\documentstyle[amssymb,amsfonts]{amsart}

 at 10 true pt

\def \endprf{\hfill  {\vrule height6pt width6pt depth0pt}\medskip}
\def\emph#1{{\it #1}}
\def\textbf#1{{\bf #1}}

\def\CS{{\mathcal S}}
\def\CT{{\mathcal T}}

\def\BBN{{\mathbb N}}

\def\BBZ {{\mathbb Z}}

\def\hs{\hfill $\square$}

\def\mod {\hbox {\rm mod}}

\def\ep{\epsilon}

\theoremstyle{plain}
  \newtheorem{theorem}[subsection]{Theorem}

  \newtheorem{fact}[subsection]{Fact}
  \newtheorem{lemma}[subsection]{Lemma}
  \newtheorem{corollary}[subsection]{Corollary}
  \newtheorem{claim}[subsection]{Claim}
  \newtheorem{defn}[subsection]{Definition}
\newtheorem{conj}[subsection]{Conjecture}
\theoremstyle{remark}

\theoremstyle{definition}

\def\ep{\epsilon}

\def\beq{\begin{equation}}
\def\eeq{\end{equation}}
\def\barray{\begin{eqnarray*}}
\def\earray{\end{eqnarray*}}

\def\CS{{\mathcal S}}
\def\CT{{\mathcal T}}

\def\BBN{{\mathbb N}}

\def\BBZ {{\mathbb Z}}
\def\BBP {{\mathbb P}}

\def\hs{\hfill $\square$}

\def\mod {\hbox {\rm mod}}
\def\Volume{\hbox {\rm Vol}}
\def \rank{\hbox{\rm rank}}
\def\Vol{\hbox {\rm Vol}}

\begin{document}

\title{ Long arithmetic progressions in sumsets: Thresholds and Bounds}

\author{E. Szemer\'edi}
\address{Department of Computer Science, Rutgers University, NJ 08854.}
\email{szemered@@cs.rutgers.edu}
\thanks{E. Szemer\'edi is supported by an NSF grant.}

\author{V.  Vu}
\address{Department of Mathematics, UCSD, La Jolla, CA 92093-0112}
\email{vanvu@@ucsd.edu}

\thanks{V. Vu is an A. Sloan  Fellow and is supported by an NSF Career Grant.}

\begin{abstract}  For a set $A$ of integers, the
sumset $lA =A+\dots+A$ consists of those numbers which can be
represented as a sum of $l$ elements of $A$

$$lA =\{a_1+\dots a_l| a_i \in A_i \}. $$

 A closely related and equally interesting
notion is that of $l^{\ast}A$, which is the collection of numbers
which can be represented as a sum of $l$ different elements of $A$

$$l^{\ast} A =\{a_1+\dots a_l| a_i \in A_i, a_i \neq a_j \}. $$

The goal of this paper is to investigate the structure of  $lA$
and $l^{\ast}A$, where $A$ is a subset of $\{1,2, \dots, n\}$. As
applications, we solve two conjectures by Erd\"os and Folkman,
posed in  sixties.
\end{abstract}

\maketitle


 \noindent Math Review Numbers: 11B25, 11P70,
11B75.


\tableofcontents

\section{  Overview}

One of the main tasks of additive number theory is to examine
structural properties of sumsets.  For a set $A$ of integers, the
sumset $lA =A+\dots+A$ consists of those numbers which can be
represented as a sum of $l$ elements of $A$

$$lA =\{a_1+\dots +a_l| a_i \in A_i \}. $$

 A closely related and equally interesting
notion is that of $l^{\ast}A$, which is the collection of numbers
which can be represented as a sum of $l$ different elements of $A$

$$l^{\ast} A =\{a_1+\dots +a_l| a_i \in A_i, a_i \neq a_j \}. $$

 Among the most well-known results in all mathematics are
Vinogradov's theorem which says that $3 \BBP$ ($\BBP$ is the set
of primes) contains all sufficiently large odd number and Waring's
conjecture (proved by Hilbert, Hardy and Littlewood, Hua, and many
others) which asserts that for any given $r$, there is a number
$l$ such that $l^{\ast} \BBN^r$ ($\BBN^r$ denotes the set of
$r^{th}$ powers) contains all sufficiently large positive integers
(see  \cite{Vau} for an excellent exposition concerning these
results).

\vskip2mm

In recent years, a considerable amount of attention has been paid
to the study of finite sumsets. Given a finite set $A$ and a
positive integer $l$, the natural analogue of Vinogadov-Waring
results is to show that under proper conditions,  the  sumset $lA$
($l^{\ast} A$) contains a long arithmetic progression.

\vskip2mm

 Let us assume  that $A$ is a subset of the
interval $[n] =\{1, \dots, n \}$, where $n$ is a large positive
integer. The concrete problem we would like to address is to
estimate the minimum length of the longest arithmetic progression
in $lA$ ($l^{\ast} A$) as a function of $l, n$ and $|A|$. We
denote this function by $f(|A|, l, n)$ ($f^{\ast}(|A|, l, n)$),
following a notation in \cite {Fresurvey}. Many estimates for
$f(|A|, l, n)$ have been discovered by Bourgain, Freiman,
Halberstam, Green, Ruzsa, and S\'ark\"ozi (see Section 2), but
most of these results focus on sets with very high density, namely
$|A|$ is close to $n$. Estimating $f^{\ast} (|A|,l,n)$ seems much
harder and not much was known prior to our study.

\vskip2mm
 In this paper, we solve both problems almost completely for a wide range of $l$ and $|A|$.
Our study reveals a surprising  fact that the functions $f(|A|,
l,n)$ and $f^{\ast}(|A|, l,n)$ are not continuous and admits a
threshold rule. We have successfully located  the threshold points
within constant errors and established the asymptotic behavior of
the functions between consecutive threshold points. It has also
turned out,  during our study, that  the sum $l^{\ast} A$ is
indeed fundamentally harder to attack than its counterpart $l A$.

\vskip2mm

Center to our study is the development of  a new, purely
combinatorial, method. This method is totally  different from
harmonic analysis methods used by most researchers and seems quite
flexible. For instance, it is easy to extend our results in many
directions. Moreover, the method carries us far beyond our
original aim of estimating lengths of arithmetic progressions,
leading to more general theorems about proper generalized
arithmetic progressions.

\vskip2mm

Our results also have some interesting applications. In
particular, we settle two forty year old conjectures of Erd\"os
\cite{erd} and Folkman \cite{folk} (respectively) concerning
infinite arithmetic progressions.

\vskip2mm

Let  us now  present  a brief introduction to the content of our
paper:

\begin{itemize}

\item  In Section 2 we present the notion of GAPs and state
Freiman's famous inverse theorem, both of which play a crucial
role in our study. In Section 3, we first describe some earlier
results on the topic. Next, we present a construction which
suggests
 a conjecture about the length of the longest arithmetic
progression in $lA$. It would be  important to keep this
construction in mind as it motivates lots of our arguments later
on.  The first main result of Section 3 confirms the conjecture
motivated by the construction. This result, among others, reveals
the surprising  fact that $f(|A|,l,n)$ is not continuous and
admits a threshold behavior. There are many threshold points and
we are able to locate them within a constant factor. The second
main result, which refines  the first one, provides  a more
general and complete picture. We can prove that $lA$ not only
contains long arithmetic progressions, but also contains large
proper generalized arithmetic progressions (a regular arithmetic
progression is a special  proper generalized arithmetic
progression of rank one; we shall use short hand GAP for
generalized arithmetic progression). In the next section, Section
4, we prove these two results. The first four subsections of
Section 4 are devoted to the development of a variety of tools,
through which we could establish a connection between our study
and inverse theorems of Freiman type. Exploiting this connection,
we complete the proofs in the final two subsections. This
concludes the first part of the paper.

\item The second part of the paper consists of two sections,
Section 5 and Section 6. In Section 5, we generalize the results
in Section 3 to sums of different sets. Instead of considering
$lA$, we consider the sum $A_1 +\dots+A_l$, where
$|A|_1=\dots=|A_l| =|A|$. Thanks to the flexibility of  our
method, we  can extend the results of Section 3 to this setting in
a relatively simple manner. Also in this part we discuss  an
application which settles a conjecture posed by Folkman in 1966.
This conjecture was considered by Erd\"os and Graham (\cite{EG},
Section 6) the most important problem in the study of subcomplete
sequences. An infinite sequence is subcomplete if its partial sums
contains an infinite arithmetic progression. Folkman conjectured
that a sufficiently dense sequence of positive integers (with
possible repetitions) is subcomplete. In Section 6, we first work
out a sufficient condition for subcompleteness and next use the
results in Section 5 to show that a sufficiently dense sequence
should satisfy  this condition.

\item Sections 7, 8 and 9 form the third part of the paper. This
part contains our strongest result whose proof is also the
 most technical. The heart of this part is Theorem 7.1, which
  extends the  results in
Section 3  to the sumset $l^{\ast} A$.   The proof comprises
several phases. In the first phase, we prove a structural property
of a set $A$ where $l^{\ast} A$ does not contain a generalized
arithmetic progression as large as we desire. This property, which
might be of independent interest, shows that such a set $A$
contains a very rigid subset which almost looks like a
 generalized arithmetic progression. The
verification of the structural lemma occupies most of Section 7.
Section 8 contains the rest of the proof, whose core consists of
an observation about proper GAPs (subsection 8.3) and a variant of
the so-called tiling technique, introduced in an earlier paper
\cite{szemvu2}. Section 9 discusses a conjecture of Erd\"os (posed
in 1962) which is related to the above mentioned conjecture of
Folkman. This conjecture was proved in an earlier paper
\cite{szemvu2} using a special case of the main result in Section
7, but here we give  a shorter proof using the general condition
worked out in Section 6. Several other applications of the main
result of this part will appear in future papers \cite{Vunew1,
Vunew2}.

\item The last part of the paper contains Section 10, in which we
extend all previous results to finite fields. We assume that $n$
is a prime and  consider arithmetic progressions modulo $n$. This
modification will lead to a natural change in the statement of the
results, but the proofs remain basically the same. We conclude
this part by mentioning an application concerning the problem of
counting zero-sum-free sets.

\end{itemize}

The paper contains several new technical ingredients, some of
which (such as the study of proper GAPs in Sections 3 and 8 and
the rank reduction argument used in Sections 3 and 4) would be of
independent interest. Our writing benefits from two earlier papers
\cite{szemvu, szemvu2}, which established several partial results
and launched the foundation of our study. Many ideas from these
two papers will be used here, frequently in more general and  more
comprehensible  forms.

\section {Inverse Theorems}

A generalized arithmetic progression (GAP) of rank $d$ is subset
$Q$ of $\BBZ$ of the following form $\{a+ \sum_{i=1}^d x_i a_i |0
\le x_i \le n_i \}$; the product $\prod_{i=1}^d n_i$ is its volume
and we denote it by $\Volume (Q)$. In fact, as two different GAPs
might represent the same set, we always consider GAPs together
with their structures. The set $(a_1, \dots, a_d)$ is called the
{\it difference} set of $Q$.

\vskip2mm

Freiman's famous inverse theorem \cite{frei} asserts that if
$|A+A| \le c|A|$, where $c$ is a constant,  then $A$ is a dense
subset of a generalized arithmetic progression of constant rank.
In fact, the statement still holds in a slightly more general
situation, when one considers $A+B$ instead of $A+A$. This was
shown by Ruzsa \cite{Ruz}, who gave a very elegant proof which was
different from  Freiman's.

\begin{theorem}  \label{freiman}
For every positive constant $c$ there is a positive integer $d$
and a positive constant $k$ such that the following holds.  If $A$
and $B$ are two subsets of $\BBZ$ with the same cardinality and
$|A+B| \le c|A|$, then $A$  is a subset of a generalized
arithmetic progression  $P$ of rank $d$ with volume at most $k
|A|$.
\end{theorem}

The most recent estimate on $k$ (as a function of $c$) is due to
Chang \cite{Chang}. In our paper, however, we shall be more
concerned with the best value of $d$ (see Lemma \ref{bilu} in
Section 3). The following result is a simple consequence of
Fremain's theorem and Pl\"uneke's theorem (for the statement of
Pl\"uneke's theorem see, e.g., \cite{Ruz}).

\begin{theorem}  \label{freiman-twoset}
For every positive constant $c$ there is a positive integer $d$
and a positive constant $k$ such that the following holds.  If $A$
and $B$ are two subsets of $\BBZ$ with the same cardinality and
$|A+B| \le c|A|$, then $A+B$  is a subset of a generalized
arithmetic progression  $P$ of rank $d$ with volume at most $k
|A|$.
\end{theorem}

For the special case when $c$ is relatively small, one can set
$d=1$. The following is a consequence of another theorem of
Freiman \cite{frei}.

\begin{lemma} \label{lev1} The following holds for all sufficiently large $m$.
  If $A$ is a set of integers
of cardinality $m$ and $|A+A| \le 2.1m $, then $A$ is a subset of
an  arithmetic progression of length $1.1m$. \end{lemma}

Again, we can replace $A+A$ by $A+B$. The following is a corollary
of a result by Lev and Smelianski (Theorem 6 of \cite{lev}).

\begin{lemma} \label{lev2} The following holds for all sufficiently large $m$.
  If $A$ and $B$ are two sets of integers
of cardinality $m$ and $|A+B| \le 2.1m $, then $A$ is a subset of
an  arithmetic progression of length $1.1m$. \end{lemma}

Both Lemmas \ref{lev1} and \ref{lev2} are relatively simple and do
not require the inverse theorem to prove.

\section {Long arithmetic progressions in $lA$}

\subsection {Some previous results}

Problems concerning arithmetic progressions in sumsets are
non-trivial and not too many results are known. In the following,
we describe some of the main results in this area. Bourgain
\cite{Bou} proved that if $|A| = \delta n$ where $\delta$ is a
positive constant then $2A$ contains an arithmetic progression of
length $e^{\epsilon {\log^{1/3} n}}$, where $\ep$ is a positive
constant depending on $\delta$.   Freiman, Halberstam and Ruzsa
\cite{FHR} consider sumsets modulo a prime and proved that

\begin{theorem}  \label{FHR} Let $n$ be a prime and $A$ a set of
residues modulo $n$, $|A| =\gamma n$, $0 < \gamma < 1$ may depend
on $n$. Let $l $ be a positive integer at least 3. Then $lA$
contains an arithmetic progression (modulo $n$) of length $\Omega
(\gamma n ^{\frac{1}{16} \gamma^{l/(l-2)} })$. \end{theorem}

 Notice that Theorem \ref{FHR} is
stated for any $\gamma$, but it is really efficient only when
$\gamma$ is relatively large. Indeed, if one wants to have $\gamma
n ^{\frac{1}{16} \gamma^{l/(l-2)} } \ge 1$ one needs to set
$$\gamma =\Omega (\frac{1}{\ln n}). $$ So Theorem \ref{FHR} does not give a
non-trivial bound in the case $|A| =o(\frac{n}{\ln n})$.
Bourgain's result and Theorem \ref{FHR} have recently been
improved by Green \cite{Gre}, but the applicable range does not
change.

 Prior to our study,
the only result (as we know of) which applies to sets with
relatively small cardinality is the following theorem, proved by
 S\'ark\"ozy \cite{sar}.

\begin{theorem}  \label{FS} There are positive  constants $c$ and $C$ such that the following holds. If
$A$ is a subset of $[n]$ and $l$ is a positive integer such that
$l|A| \ge Cn$, then $lA$ contains an arithmetic progression of
length $cl|A|$. \end{theorem}

\noindent Answering a question of S\'ark\"ozy, Lev \cite{levsar}
shown that one can set $C$ equal to 2, which is the optimal value.

\vskip2mm

It is clear that Theorem \ref{FS} is sharp, up to a constant
factor.  Let $A$ be the set of all positive integers from $1$ to
$|A|$. Then $lA$ is the set of all positive integers from $l$ to
$l|A|$.

\vskip2mm The main result of this section gives a sharp estimate
for a wide range of $|A|$ and $l$, including  Theorem \ref{FS} as
a special case. More importantly, our proof reveals the structures
of those sets $A$ whose sumsets $lA$ do not contain a very long
arithmetic progression. In the next subsection, we describe the
construction that  motivates our result.

\vskip2mm To conclude this subsection, let us mention that the
proofs of all results mentioned in this paper, with the exception
of S\'ark\"ozi's proof, are analytic, making heavy use of harmonic
analysis,  and are very different from the proofs in this paper.

\subsection { Sudden jumps}

Our first crucial observation is that the statement of Theorem
\ref{FS} stops to hold when $l|A|$ becomes a little bit less than
$ n$. The following construction shows that there is a set $A
\subset [n]$ and a number $l$ such that $l|A| \approx n/4$ while
the length of the longest arithmetic progression in $lA$ is only
$O(l|A|^{1/2})$ (here and later $\approx$ means ``approximately").

\vskip2mm

{\it \noindent  The construction.} Let $A=\{ p_1 x_1 + p_2 x_2| 1
\le x_1 \le m \}$, where $p_1 \approx p_2 \approx \frac{n}{2m} $
are two primes
 and $p_2 > m$.
 It is convenient to think of $A$ as a square in the two dimensional lattice
 $\BBZ^2$. A point $(x_1,x_2)$ corresponds to the number
$p_1x_1+p_2x_2$. It is easy to show that this correspondence is
one to one. Indeed

$$p_1x_1+p_2x_2 = p_1x_1'+p_2x_2' $$

\noindent implies that

$$p_1(x_1-x_1') = p_2(x_2'-x_2) $$

\noindent which is impossible because of divisibility and the fact
that $|x_1-x_1'| < m < p_2$. Thus, $|A|= m^2$. Let $l=
\frac{n}{(4+\ep)|A|} = \frac{n}{(4+\ep)m^2}$, where $\ep$ is an
arbitrary positive constant. We have

 $$lA= \{ p_1 x_1 + p_2 x_2| l \le x_1 \le lm \}. $$

Let $P$ be an AP in $lA$, we are going to show that the
coordinates of the elements of $P$ also form an AP of the same
length. Thus $|P|$ is at most the length of an edge of $l|A|$,
which is less than $lm= l|A|^{1/2}$. Observe that

\centerline { $p_2 \approx n/2m \ge 2lm  \,\, \hbox{ since}\,\, l=
n/(4+\ep)|A|= n/(4+\ep)m^2. $}

Consider three consecutive terms in $P$, $z+z^{''} = 2z'$. Write
$z= p_1x_1+ p_2 x_2$. We have

$$  (p_1x_1+ p_2 x_2)  +   (p_1x_1^{''} + p_2 x_2^{''} ) = 2 (p_1x_1' + p_2 x_2'), $$

\noindent which implies

$$ p_1 (x_1 + x_1^{''} - 2x_1')  = - p_2  (x_2 + x_2^{''} - 2x_2'), $$

\noindent which is again impossible as

$$ |x_1 + x_1^{''} - 2x_1'| < 2lm  \le p_2. $$

\vskip2mm

\noindent Next, we generalize the above construction to higher
dimensions.

\vskip2mm

{\it \noindent  The general  construction.} Let $d$ be a constant
positive integer at least 2 and $\delta$  be a small positive
constant. Consider two numbers $|A|$ and $l$ satisfying
$l^{d-1}|A| \le \frac{1 -\delta}{2d} n$. We shall construct a set
$A$ of cardinality $|A|$ such that the longest arithmetic
progression in $lA$ has length $l|A|^{1/d}$.

\vskip2mm  Set $a= \lfloor \frac{(1-\delta/3)n}{d|A|^{1/d} }
\rfloor$ and $b = \lfloor (\frac{n}{d |A|^{1/d}})^{1/(d-1)}
\rfloor $. Set $b_1=0, b_2=1$ and if $d \ge 3$ then set $b_i=
\lfloor b ^{(i-2)/(d-1)} \rfloor$ for all $3 \le i \le d$. Finally
set $a_i=a +b_i$. It is a routine to verify that for a
sufficiently large $n$

\begin{equation}  \label{const1}  (1-\delta/3) a^{1/(d-1)} \ge 2l|A|^{1/d}.
\end{equation}

Consider the set

$$A = \{\sum_{i=1}^d a_i x_i| 1 \le x_i \le
|A|^{1/d}\} $$ (for convenience we assume that $|A|^{1/d}$ is an
integer). The term $\frac{1-\delta}{d}$ in the definition of $a$
guarantees that $P$ is a subset of $[n]$. It is convenient to view
both $A$ and  of $lA$ as $d$-dimensional integral boxes. The edges
of $lA$ form arithmetic progressions of length $l |A|^{1/d}$.
Similar to the case $d=2$, we are going  to prove the following
two claims.

\begin{claim} \label{nolonger} $lA$ does not contain an arithmetic progression of
length larger than $l |A|^{1/d}$. \end {claim}

\begin {claim} \label{card} The cardinality of $A$ is $|A|$.
\end{claim}

\vskip2mm

{\bf \noindent Proof of Claim \ref{nolonger}.}  Consider  an
arithmetic progression $P$ in $lA$ and let $z, z^{'}, z^{''}$ be
three consecutive elements of $P'$. We have $z +z ^{''}=2z^{'}$.
Write $z= \sum_{i=1}^d a_i x_i, z'= \sum_{i=1}^d a_i x^{'}_i$ and
$z^{''}=\sum_{i=1}^d a_ix^{''}_i$, it follows that $\sum_{i=1}^d
(x_i+x^{''}_i-2x^{'}_i) a_i =0$. Notice that $1 \le x_i, x^{'}_i,
x^{''}_i \le l|A|^{1/d}$, so  $|x_i +x^{''}_i -2x^{'}_i| < 2
l|A|^{1/d}$, for all $i$'s.

\vskip2mm

Next, we  show that the diophantine equation $\sum_{i=1}^d r_i a_i
=0$ cannot have non-trivial roots with small absolute values,
namely, $|r_i| < 2l|A|^{1/d}$ cannot hold simultaneously for all
$i$'s. Consider a non-trivial root $\{r_1, \dots, r_d \}$. There
are two cases

\vskip2mm

(I) $\sum_{i=1}^d r_i =0$. By the definition of the $a_i$'s, it
follows that $\sum_{i=1}^d r_i b_i=0$ and $d$ should be at least
3. Let $j$ be the largest index where $r_j \neq 0$, it is easy to
see that $j \ge 3$.  On the other hand, by the definition of the
$b_i$'s,  for any $j \ge 3$

\begin{equation}  \label{lower1} \max_{ 1\le i \le d} |r_i| \ge \frac{b_j}
{\sum_{i=1}^{j-1} b_i }\ge a^{1/(d-1)} \ge 2l|A|^{1/d},
\end{equation}

\noindent where the last inequality is from (\ref{const1}).
\vskip2mm

(II) $\sum_{i=1}^d r_i \neq 0$. In this case, it is obvious that

\begin{equation}  \label{lower2} \max_{1\le i \le d} |r_i| \ge  \frac{a}
{\sum_{i=1}^{d} b_i} \ge  (1-\delta) a^{1/(d-1)} \ge 2l|A|^{1/d}.
\end{equation}

\vskip 2mm

By the previous facts, we can conclude that
$x_i+x^{''}_i-2x^{'}_i=0$ for all $i$'s. So for each $i$, the
coordinates of $z_i$ form an arithmetic progression. This implies
that the length of $P$ could be at most the length of the ``edges"
of $A$, which is $l|A|^{1/d}$. \hs

\vskip2mm

From the previous proof, it is obvious that if $\sum_{i=1}^d a_i
x_i = \sum_{i=1}^d a_i x'_i $ for $ 1 \le x_i, x'_i \le |A|^{1/d}$
for all $1 \le i \le d$, then $x_i =x'_i$ for all $i$'s. This
implies that  the cardinality of $A$ is $|A|$, proving Claim
\ref{card}.

\vskip2mm

This construction plays a very important role in the whole paper.
It not only leads us to the statements of our theorems, but also
motivates many of our arguments.

\vskip2mm

{\it \noindent  The sudden jumps.} For the sake of simplicity, let
us consider $l$ and $n$ fixed and view $f(|A|, l,n)$ as a function
of $|A|$ (we call this function $g(|A|)$). The special case $d=2$
shows that if $|A| \le \frac{1-\delta}{4} \frac{n}{l}$, then
$g(|A|)$ is upper bounded by $l|A| ^{1/2} $. This and Theorem
\ref{FS} imply that $g(|A|)$ admits a dramatical  change in order
of magnitude somewhere near the point $\frac{n}{l}$. If $|A| \ge C
\frac{n}{l}$ for some sufficiently large constant $C$, then
$g(|A|)$ (up to a multiplicative constant) behaves like $l|A|$. On
the other hand, if $|A| \le \frac{1-\delta}{4} \frac{n}{l}$ then
then $g(|A|)$ is upper bounded by $l|A| ^{1/2} $. This indicates
that $g(|A|)$ is not a continuous function and its behavior must
follow a threshold rule.

\vskip2mm

The general construction  suggests that  $n/l$ is not the only
threshold (a place where $g(|A|)$ jumps). Assume, for a moment,
that we could prove that close to the left of $n/l$, $g(|A|)$
behaves like $l|A|^{1/2}$. This behavior, however,  cannot
continue to hold with $|A|$ getting significantly smaller than
$n/l$. Indeed, once $|A| $ becomes less than  $\frac{1-\delta}{6}
\frac{n}{l^2}$ then $g(|A|)$ is upper bounded by $l |A|^{1/3}$.
Thus, another threshold should  occur around the point $
\frac{n}{l^2}$. Motivated by this reasoning, one would conjecture
that there is a threshold around $\frac{n}{l^d}$ for any fixed
positive integer $d$. To the right of the threshold, $g(|A|)$
behaves like $l |A|^{1/d}$, while to the left  it  behaves like $l
|A|^{1/(d+1)}$.

\subsection {$g(|A|)$ must jump}

\noindent Our first main result  confirms the above conjecture.

\begin{theorem}  \label{1} For any fixed positive integer $d$
there are positive constants $C$ and $c$ depending on $d$
 such that the following holds. For any positive integers $n $
and $l$ and any set  $A \subset [n]$ satisfying $l^{d} |A| \ge C
n$,  $lA$ contains an arithmetic progression of length $cl
|A|^{1/d}$.
\end{theorem}



\begin{corollary} \label{11} For any fixed positive integer $d$  there are positive
constants $C_1, C_2$, $c_1$ and $c_2$ depending on $d$ and $\ep$
such that whenever $\frac{C_1n}{l^{d}} \le |A| \le \frac{C_2
n}{l^{d-1}}$

$$c_1 l |A|^{1/d} \le  f(|A|, l, n) \le c_2 l |A|^{1/d}. $$

\end{corollary}

\noindent  Let us again consider $f(|A|,l,n)$ as a function
$g(|A|)$ of $|A|$,
  assuming $n$ and $l$ are fixed. It is more convenient to view $g(|A|)$
  on a
logarithmic scale. For this purpose, let us define $x= \ln |A|$
and $y (x) = \ln g(|A|)$. Corollary \ref{11} implies

\begin{corollary} \label{11-0} For any fixed positive integer $d$  there are
constants $C_1, C_2$, $c_1$ and $c_2$ depending on $d$ such that
whenever $ \ln n -d \ln l +C_1 \le x \le  \ln n -(d-1) \ln l +C_2$

$$  \frac{1}{d}  x + \ln l + c_1 \le  y(x)  \le \frac{1}{d}  x  + \ln l + c_2. $$

\end{corollary}

\noindent The values of the constants $C_1, C_2, c_1, c_2$ in this
corollary are, of course, different from the values of $C_1, C_2,
c_1, c_2$ in Theorem \ref{1}. Corollary \ref{11-0} determines the
value of $y(x)$ up to a constant additive term for all $x$ except
few intervals of constant lengths. An exceptional interval is a
neighborhood of  a threshold point $\ln n - d\ln l = \ln
\frac{n}{l^d}$ and is of the form $[\ln n - d \ln l +C_2(d-1), \ln
n - d \ln l+C_1(d)]$, which has length $C_1(d)-C_2 (d-1)$. Here we
write $C_1(d)$ and $C_2(d-1)$  instead of $C_1$ and $C_2$ to
emphasize the dependence on $d$ and $d-1$, respectively.

\vskip2mm

The above results  locate the thresholds within constant factors.
It would be nice to find the exact locations of these thresholds.

\vskip2mm

{\bf \noindent Question.} {\it Find the exact values of the
constants $C$ and $c$ in Theorem \ref{1}.} \vskip2mm

\noindent The case $d=1$ was treated by Lev in \cite{levsar}. For
general $d$, our  construction shows that $C(d)$ is at least
$(1-o(1))/2d$.

\subsection {A stronger theorem about generalized arithmetic progressions }

Theorem \ref{1} is a only a tip of an iceberg and we are going to
extend it  in various directions. In the first extension, we show
that Theorem \ref{1} is a consequence of a  stronger theorem about
GAPs.

\vskip2mm In order to guess what we may say about the possible
existence of GAPs in $lA$, let us go back to the construction.
Observe that the constructed sumset $lA$ contains not only an
arithmetic progression of length $l|A|^{1/d}$, but  also  a proper
GAP of rank $d$ and cardinality $\Omega (l^d |A|)$. The arithmetic
progression of length $l|A|^{1/d}$ we talked about is actually an
edge of this  GAP. Thus, our first guess is, naturally, that $lA$
contains a GAP of rank $d$ and cardinality $\Omega (l^d |A|)$.
This guess is, nevertheless,   false. To see this, notice  that if
we let $A$ in the construction be a GAP of dimension $d' < d$ with
appropriate parameters, then $lA$ is a GAP of dimension $d'$ of
cardinality $\Omega (l^{d'} |A|) $ which is much less than $\Omega
(l^d |A|)$ (it is interesting to note that  in this case $lA$
contains an arithmetic progression of length $\Omega (l|A|^{1/d'})
\gg \Omega (l |A|^{1/d})$). So, the strongest statement one could
say is that $lA$ contains a proper GAP of rank $d'$ and
cardinality $\Omega (l^{d'} |A|)$ for some integer $1 \le d' \le
d$. This  turns out to be the truth.

\begin{theorem}  \label{111} For any fixed positive integer $d$
there are positive constants $C$ and $c$ depending on $d$
 such that the following holds. For any positive integers $n $
and $l$ and any set  $A \subset [n]$ satisfying $l^{d} |A| \ge C
n$,  $lA$ contains a proper GAP of rank $d'$ and volume at least
$cl^{d'} |A|$, for some integer $1 \le d' \le d$.
\end{theorem}

The other main results of this paper, Theorems \ref{211},
\ref{311}, \ref{411}, \ref{211modn} are extensions of this theorem
in various directions.

\vskip2mm  To conclude this subsection, let us point out that both
Theorem \ref{1} and Theorem \ref{111} are invariant under affine
transformations. Instead of assuming that $A$ is a subset of
$[n]$, we can assume that $A$ is a subset of an arithmetic
progression of length $n$. In fact, for technical reasons, we will
frequently assume that $A$ contains $0$.

\subsection {More about generalized arithmetic progressions}

\vskip2mm  Consider a GAP $Q= \{a+ \sum_{i=1}^d x_ia_i| 0 \le x_i
\le n_i\}$. It is convenient to consider $Q$ together with the box
$B_Q =\{(x_1, \dots, x_d) | \ 0 \le x_i \le n_i\}$  of $d$
dimensional vectors and the following map $\Phi$ from $\BBZ^d$ to
$\BBZ$

$$\Phi (x_1, \dots, x_d) = a+ \sum_{i=1}^d x_i a_i. $$

\noindent The volume of $Q$ is the geometrical volume of the
$d$-dimensional box spanned by $B_Q$

$$\Volume (Q) =\Volume
(B_Q) = \prod_{i=1}^d n_i. $$

 We say that $Q$ is {\it proper} if $\Phi (B_Q)$ is
injective. In this case the cardinality of $Q$ is $\prod_{i=1}^d
(n_i+1) = |B_Q|$. It is trivial that

\begin{equation}  \label{gap1} |Q| \le 2^d \Volume (Q), \end{equation}

\noindent and if $Q$ is proper then

\begin{equation}  \label{gap2} \Volume (B_Q) < |B_Q| \le 2^d \Volume (B_Q).
\end{equation}

\vskip2mm If $Q$ is not proper, then there are two vectors $u$ and
$w$ in $B_Q$ such that $\Phi (u) =\Phi(w)$. The vector $v=u-w$ is
called a {\it vanishing} vector. By linearity, it is clear that if
$v$ is vanishing then $\Phi (v)=0$ and $\Phi (v+ u) =\Phi (u)$ for
any $u \in \BBZ^d$.

\vskip2mm

In the following we specify some rules used in calculation
involving GAPs.

\vskip2mm

{\it \noindent Addition.} We only add two GAPs with the same
difference set and the result is a GAP with this difference set.
For instance, if $P= \{ a+ a_1x_1+\dots a_d x_d | 0 \le x_i \le
m_i\}$ and $Q= \{ b+ a_1x_1+\dots a_d x_d | 0 \le x_i \le n_i\}$
then

$$P+Q =\{(a+b) + a_1x_1+\dots a_d x_d | 0 \le x_i \le m_i+n_i\}. $$

\noindent Substraction is defined similarly.

\vskip2mm

{\it \noindent Multiplication.} For a GAP $P$, we have $2P = P+P$
and $lP = (l-1)P+ P$.  \vskip2mm

{\it \noindent Division.} Consider a GAP  $P= \{ a+ a_1x_1+\dots
a_d x_d | 0 \le x_i \le m_i\}$. We say $P$ is {\it normal} if
$a=0$. In this case, we define

$$\frac{1}{s} P = \{  a_1x_1+\dots
a_d x_d | 0 \le x_i \le m_i/s\}. $$

\noindent All of our arguments concerning GAPs are  invariant with
respect to affine transformation (shiftings in particular), so we
could (and shall) automatically assume that a GAP is normal when
it is involved in  division.


\subsection {Some simple tricks}

In this subsection, we describe several simple tricks which we use
frequently throughout the paper.

\vskip2mm

As $C$ can be set arbitrary large, we can sacrify constant factors
in many arguments. So we are going to make several assumptions,
whose ``prices" are only constant factors,   which are very
convenient for the proofs.

\vskip2mm

{\it \noindent Divisibility.}  By increasing the value of $C$, we
can assume that $l$ is a power of two. Indeed, if we replace $l$
by the closest power of two, then the magnitude of $l$ decreases
by at most 2. Similarly, once we have  a GAP of constant rank and
all we care is the volume of this GAP, up to a constant factor,
then we can assume that the lengths of the edges are divisible by
2 (or by any fixed integer).
 This latter
assumption is convenient for divisions. For instance, whenever we
need to divide a GAP $P$ by a constant $s$, we shall always assume
that the lengths of the edges of $P$ are divisible by $s$.

\vskip2mm

{\it \noindent Passing to subsets.} In many situations, it is
useful to assume  that  a certain set, say $X$, has a certain
property. On the other hand,  we can only prove that $X$ has a
subset $X'$ with the desired property. However, when $X'$ has
constant density in $X$, we can frequently assume that $X$ has the
desired property, again by increasing the value of $C$.

\vskip2mm

{\it \noindent A graph with small degrees contains a large
independent set.} A graph consists of a set $V$ of vertices and a
set $E$ of edges, where an edge is a pair of two different
vertices. The degree of a vertex $v$ is the number of edges
containing $v$. If $(u,v)$ is an edge, then $u$ is a neighbor of
$v$ and vice versa. A subset of $V$ is called {\it independent} if
it does not contain any edge. We are going to use the following
simple fact from graph theory.

\begin{fact} \label{independentset} Let $G$ be a graph on $n$
vertices. Assume that any vertex of $G$ has degree at most $d$.
Then $G$ contains an independent set of size $n/(d+1)$. \end{fact}

{\bf \noindent Proof.} Let $I$ be a maximal independent set. Since
$I$ is maximal, the neighbors of the vertices in $I$ and $I$
together cover the vertex set of $G$. Since the vertices of $I$
have at most $d|I|$ neighbors, it follows that

$$d|I| +|I| \ge n, $$

\noindent proving the claim. \hs

\vskip2mm The above fact implies that if $G$ does not contain an
independent set of size $s$, then $G$ has a vertex with degree at
least $n/s$.

\section {Proofs of Theorem \ref{1} and Theorem \ref{111}}

This section has six subsections. In the first four subsections we
develop a variety of tools. The proof of Theorem \ref{1} and that
of Theorem \ref{111} are presented in the last two subsections.

\vskip2mm

Let us start with a sketch of the proof of Theorem \ref{1}.
Consider the sequence $$A, 2A, 4A,..., lA$$(without loss of
generality we can assume that $l$ is a power of 2). Since $lA$ is
a subset of the interval $[ln]$, $|lA|$ is at most $ln$. This
implies that the ratio $|2^{i+1} A|/|2^i A|$ cannot always be
large.  In particular, there is a constant  $K$ such that
$|2^{i+1} A| \le K |2^i A|$ holds for some index $i $ less than
$\log_2 l$. On the other hand, $2^{i+1} A = 2^iA + 2^iA$, so by
applying Freiman's theorem we can deduce that $2^iA$ is a dense
subset of a GAP $P$ with constant rank.

\vskip2mm

Let us assume, for a moment,  that $2^iA$ has density one in $P$,
namely, $2^iA=P$. Thus $2^iA$ contains a long arithmetic
progression $B$ of length at least $(\Vol P)^{1/rank(P)}$. As $i$
is  less than $\log_2 l$, $lA$ contains an even longer arithmetic
progression of length at least $\frac{l}{2^i} |B|$.

\vskip2mm

In order to carry out this scheme, we first need to show that
assuming $2^iA =P$ is not oversimplifying. This will be carried
out in the second subsection, where we show that at the cost of
constant factors we can think of a dense subset of a GAP as the
whole set.

\vskip2mm

With the aid of this assertion, it is now not so hard to prove
that $lA$ contains an arithmetic progression of length $l
|A|^{\ep}$ for some small $\ep$. In order to optimize $\ep$, we
need to optimize $K$ and the rank of $P$. The optimal value of $K$
is easy to guess while the optimal value of the rank of $P$ will
be provided by a result of Bilu \cite{bil}, which is a part of his
proof of Freiman's theorem.

\vskip2mm

Now comes the last, and perhaps most intriguing point. Even with
these optimal parameters, we could not obtain the bound claimed in
the theorem (however, we can obtain a weaker theorem proved in an
earlier paper \cite{szemvu}). To fill in the gap, we need to prove
certain properties of non-proper and proper GAPs. These properties
lead us to Lemma \ref {t1-7} which is the main lemma of the proof.
The verification of this lemma requires the preparation carried
out throughout the first three subsections.

\vskip2mm

Now let us say something about the proof of Theorem \ref{111}. The
first step  is to realize that we can assume that $2^iA$ is not
only a GAP, but also a proper one. The sumset $lA$ contains a
multiple of this GAP. The trouble is that a multiple of a proper
GAP does not need to be  proper. What saves us here is a technique
called "rank reduction". The heart of this technique is an
argument which shows that under certain circumstances a multiple
of a proper GAP either is proper or contains a proper GAP of
strictly smaller rank and comparable cardinality. Thus if we fail
to complete our task in the first attempt, we can pass to a proper
GAP with smaller rank and make a new try. The GAP we start with
has a constant rank so sooner or later we must be done. The reader
would notice that this approach, in spirit, is consistent with the
statement of Theorem \ref{111}, which confirms the existence of a
GAP of rank $d'$ where $d'$ is an undetermined quantity between 1
and $d$. This value $d'$ is exactly  where the rank reduction
terminates.

\subsection {A property  of non-proper GAPs}

\vskip2mm Let us consider the ratio between the cardinality and
the volume of a GAP $P$. Assume that $P$ has the form $P= \{a +
a_1 x_1 + \dots a_d x_d| 0 \le x_i \le n_i\}$, where all $n_i's
\ge 1$. The volume of $P$ is $\prod_{i=1}^d n_i$. If $P$ is
proper, then its cardinality is $\prod_{i=1}^d (n_i +1)$ and the
ratio in question is $\prod_{i=1}^d (1+ \frac{1}{n_i})$, which is
a number between 1 and  $2^d$. For a non-proper GAP,  it is safe
to
 say  that the ratio is  less than $2^d$, but it could still be larger than 1.
  We are going to
show, nevertheless, that if $P$ is a sufficiently large multiple
of a non-proper GAP, then this ratio is bounded from above by any
fixed positive constant $\ep$.

\begin{lemma} \label{t1-4} For any positive constants
$\ep$ and $d$ there is a constant $g$ such that  the following
holds. If a GAP $Q$ of rank $d$ is not proper, then $|gQ| \le \ep
\Volume (gQ)$. Moreover,

$$|2Q| \le (1-\frac{1}{2^{d+1}}) |2B_{Q}|.$$

\end{lemma}


In the proof, we are going to use  terminologies introduced in
subsection 3.13. The reader may want to read this subsection again
before checking the proof.

\vskip2mm

{\bf \noindent Proof of Lemma \ref{t1-4}.} We can assume that $Q
=\{x_1a_1 + \dots + x_d a_d| 0 \le x_i \le n_i\}$. We consider $Q$
together with the box $B_Q$ and the canonical  map $\Phi$ from
$B_Q$ to $Q$. Since $Q$ is not proper, there is a vanishing vector
$v$ where $ -n_i \le v_i \le n_i$ for all $i=1, \dots, d$. Without
loss of generality, we can assume that the first $d'$ coordinates
of $v$ is positive and the remaining ones are non-positive. Thus
$0 < v_i \le n_i$ for $i=1, \dots, d'$ and $-n_i \le v_i \le 0$
for $d'< i \le d$. \vskip2mm

Let $h < g$ be sufficiently large integers and let  $B'$ be the
set of vectors $w$ in $gB_Q$ such that $w+v, w+ 2v, \dots, w+ hv$
are also in $gB_Q$. As $v$ is vanishing $\Phi (w) = \Phi (w+ v)
=\dots = \Phi (w+ hv)$. It follows that

\begin{equation}  \label{BQ} |gQ| \le |gB_Q \backslash B'| + \frac{1}{h+1}
|B'|= |gB_Q| -\frac{h}{h+1} |B'|,\end{equation}



\noindent  which implies

\begin{equation}  \label{BQ0} |gQ| \le (1- \frac{h}{h+1} \frac{|B'|}{|gB_Q|})
|gB_Q| \le  2^d (1- \frac{h}{h+1} \frac{|B'|}{|gB_Q|}) \Volume
(gB_Q),\end{equation}

\noindent where in the last inequality we use the trivial fact
that $|gB_Q| \le 2^d \Volume (gB_Q)$ (see (\ref{gap1}).  Next we
bound $|B'|$ from below. A vector $w$ is surely in $B'$ if $0 \le
w_i \le (g-h) n_i$ for $i \le d'$ and $hn_i \le w_i \le g n_i$ for
$ d'< i\le d$. Thus the cardinality of $B'$ is at least
$\prod_{i=1}^d \big((g-h)n_i + 1\big) $. Moreover, $|gB_Q| \le
\prod_{i=1}^d (gn_i +1)$, so

\begin{equation}  \label{BQ2} \frac{h}{h+1} \frac{|B'|}{|gB_Q|} \ge
\frac{h}{h+1} \prod_{i=1}^n \frac{ \big((g-h)n_i + 1\big)}{gn_i
+1}. \end{equation}

\noindent For any given $\ep,d$  we could choose $g$ and $h$
(depending only $\ep$ and $d$) so that

$$\frac{h}{h+1} \prod_{i=1}^d \frac{ \big((g-h)n_i + 1\big)}{gn_i
+1} \ge 1 -\ep/2^d, $$

\noindent holds for any positive integers $n_i$'s. With this
choice of $g$ and $h$, the right most formula  in (\ref{BQ0}) is
thus at most $\ep \Volume (gQ)$, proving the first statement of
the lemma. To verify the second statement, set $g=2$ and $h=1$. We
obtain

\begin{equation}  |2Q| \le (1- \frac{1}{2} \frac{|B'|}{|gB_Q|}) |gB_Q| \le
(1-\frac{1}{2} \prod_{i=1}^d \frac{n_i+1}{2n_i +1})  |2B_Q|.
\end{equation}

\noindent The product  $\prod_{i=1}^d \frac{n_i+1}{2n_i +1}$ is
larger than $\frac{1}{2^d}$ so it follows that

\begin{equation}  |2Q| \le (1- \frac{1}{2^{d+1}}) |2B_Q|, \end{equation}

\noindent completing the proof.  \hs

\subsection {The proper filling lemma}

\noindent In this subsection, we present several  lemmas which
allow us to think of a dense subset of a GAP as the whole set, at
the cost of constant factors. The first such lemma was proved in
\cite{szemvu}.

\begin{lemma} \label{full10}  For any positive constant $\gamma$ and
any positive integer $d$ there is a constant positive integer $h$
and a positive constant $\gamma'$  depending on $\gamma $ and $d$
such that the following holds. If $P$ is a generalized arithmetic
progression of rank $d$ and $B$ is a subset of $P$ such that $|B|
\ge \gamma \hbox {\rm Vol} (P)$, then $hB$ contains a generalized
arithmetic progression of rank $d$ with cardinality at least
$\gamma'|B|$.
\end{lemma}

\noindent We call this lemma the ``filling lemma", as our
motivation is to fill out a complete GAP.  Next, we strengthen
this lemma by adding a requirement that the GAP contained in $hB$
must be proper.

\begin{lemma} \label{full1}  For any positive constant $\gamma$ and
any positive integer $d$ there is a constant positive integer $h$
and a positive constant $\gamma'$  depending on $\gamma $ and $d$
such that the following holds. If $P$ is a generalized arithmetic
progression of rank $d$ and $B$ is a subset of $P$ such that $|B|
\ge \gamma \hbox {\rm Vol} (P)$, then $hB$ contains a proper
generalized arithmetic progression of rank $d$ with cardinality at
least $\gamma'|B|$.
\end{lemma}


We shall, naturally, refer to Lemma \ref{full1} as the ``proper
filling lemma''. The proof of Lemma \ref{full1} combines Lemma
\ref{full10} with the result of the previous subsection.

\vskip2mm

{\bf \noindent Proof of Lemma \ref{full1}.} By Lemma \ref{full10},
$hB$ contains a GAP $Q$ with cardinality $\Omega (|B|)$. It
suffices to  show that $Q$ contains a proper GAP  of the same rank
with cardinality $\Omega (|Q|)$. As $h=O(1)$, $\Volume (hP) = O(
\Volume (P)) = O(|B|) $, so we can assume that

\begin{equation}  \label {full1-1} |Q| \ge \gamma_1 \Volume (hP) \end{equation}

\noindent  for some positive constant $\gamma_1$.

\vskip2mm   Let $g$ be a large constant integer. Without loss of
generality we can assume that $Q =\{x_1a_1 +\dots x_d a_d | 0 \le
x_i \le n_i \}$ and $n_i$ is divisible by $g$. Let $\ep$ be a
positive constant smaller than $\gamma_1$ and consider the GAP
$Q'= \frac{1}{g} Q$. If $Q'$ is proper then we are done as

$$|Q'| \ge \Volume (Q') = \Omega (\Volume( Q)) = \Omega (|Q|). $$

\noindent  We next show that $Q'$ is indeed proper given that $g$
is sufficiently large. Assume otherwise. Choosing $g$ as in Lemma
\ref{t1-4} we have

\begin{equation}   |Q| = |gQ'|  \le \ep \Volume (gQ') = \ep \Volume (Q) \le \ep
\Volume (hP) < \gamma_1  \Volume (hP), \end{equation}

\noindent which contradicts (\ref{full1-1}). This completes the
proof. \hs

\subsection {$(\delta,d)$-sets}

We begin this subsection with an important definition.

\begin{defn} \label{t1-1} A set $A$ is a $(\delta, d)$-set if
one can find a GAP $Q$ of rank $d$ such that $B= Q \cap A$
satisfies $|B| \ge \delta \max \{|A| ,\Volume (Q)\}$. \end{defn}

\noindent The filling lemmas tell us that  a $(\delta, d)$-set
(where both $\delta$ and $d$ are constant) can be  treated as a
GAP of rank $d$, if we are allowed to sacrifice constant factors.

\begin{lemma} \label{t1-2} For any positive constants $\delta$ and
$d$ there are positive constants $g$ and $\gamma$  such that the
following holds. If $A$ is a $(\delta, d)$-set then $gA$ contains
a proper GAP of rank $d$ with cardinality at least $\gamma|A|$.
\end{lemma}

\noindent Now we are going to present another  lemma, which
supplies a sufficient condition for a set to be a
$(\delta,d)$-set. In order to motivate  this lemma, let us go back
to Freiman's inverse theorem. Freiman's theorem shows that if
$|A+A| \le c|A|$, then $A$ is a dense subset of a GAP $P$ of rank
$d = d(c)$. As we mentioned at the beginning of this section, the
optimal value of $d$ is critical to us. Observe that if $A$ is a
proper GAP of dimension $d$, then $|A+A| \le 2^d |A|$. So, one may
wonder  whether one can set $d = \lfloor \log_2 c \rfloor$.
Unfortunately, Freiman's theorem is not true with this value of
$d$ (the best known bound is $d= \lfloor c \rfloor$). On the other
hand, if we can afford  to sacrifice constant factors, then we can
actually obtain this optimal value of $d$. To be more precise, if
$|A+A| \le c |A|$, then  a constant fraction of $A$ is contained
in a GAP $P$ of ranked $d = \lfloor \log_2 c \rfloor$ with small
volume. The following lemma is a consequence of Theorem 1.3 of
\cite{bil}.

\begin{lemma} \label{bilu} For any positive constants $\ep$ and
$d$ there is a positive constant $\delta$ such that the following
holds. If $|A+A| \le (2^d-\ep)|A|$ then $A$ is a $(\delta,
d)$-set.
\end{lemma}

This lemma is a co-product of the proof of Freiman's theorem given
by Bilu in \cite{bil}.

\subsection {Rank reduction}

Now we are in position to develop the so-called rank reduction
technique, mentioned earlier in the beginning of this section.
This technique plays an important role not only in the proofs of
Theorems \ref{1} and \ref{111}, but also in the proof of Theorem
\ref{311}.

\vskip2mm

The rank reduction technique allows us to pass from one GAP to
another which has   strictly smaller rank and comparable
cardinality. We are going to present  several lemmas which
constitute the technique. The first lemmas is as follows.

\begin{lemma} \label{t1-5} For any positive constant
$d$ there is a positive constant $\delta$ such that the following
holds. If a  GAP $Q$ of rank $d$ is proper but $2Q$ is not, then
$2Q$ is a $(\delta, d-1)$-set.
\end{lemma}

{\bf \noindent Proof of Lemma \ref{t1-5}.} Applying the second
statement of  Lemma \ref{t1-4} to $2Q$ we have that

\begin{equation}  |4Q|= |2(2Q)| \le  (1-\frac{1}{2^{d+1}}) |4B_Q| \le
(1-\frac{1}{2^{d+1}}) 4^d |B_Q|, \end{equation}

\noindent where in the last inequality we used the fact that
$|4B_Q| \le 4^d |B_Q| $. Since as $Q$ is proper $|B_Q|= |Q|$.  It
follows that

\begin{equation}  |4Q| <  (1-\frac{1}{2^{d+1}}) 4^d |Q| \le (2^d -\gamma)^2
|Q|,\end{equation}

\noindent for some constant $\gamma=\gamma (d)$. It follows that
either $|2Q| \le (2^d-\gamma) |Q|$ or $|4Q| \le (2^d-\gamma)
|2Q|$. In the first case $Q$ is a  $(\delta, d-1)$-set; in the
second case $2Q$ is a  $(\delta, d-1)$-set (both statements follow
immediately from Lemma \ref{bilu}). But $Q$ is a translation of a
subset of $2Q$, so in both cases $2Q$ is a $(\delta, d-1)$-set
(notice that the three $\delta$'s in the last two sentences might
have different values). \hs

\vskip2mm

\noindent  The previous lemma and Lemma \ref{t1-2} together yield

\begin{lemma} \label{t1-6} For any positive constant
$d$ there are positive constants $g$ and $\gamma$ such that the
following holds. If a GAP $Q$ of rank $d$ is  proper but $2Q$ is
not, then $gQ$ contains a proper GAP of rank $(d-1)$ with
cardinality at least $\gamma |Q|$.
\end{lemma}

\noindent We are now ready to present the main lemma of the proofs
of Theorems \ref{1} and \ref{111}.

\begin{lemma} \label{t1-7} For any positive constants
$\ep$ and $d$ there are positive constants $c$ and $\gamma$ such
that the following holds. Let $Q$ be a proper GAP of rank $d$ and
assume that there are positive integers $l_1=2^{s_1}$ and $m$
satisfying $l_1Q \subset [m]$ and $l_1^d |Q| \ge c m$. Then there
is a positive integer $l_1'= 2^{s_1'} < l_1$ such that $l_1'Q$
contains a proper GAP $Q'$ of rank $(d-1)$ where $|Q'| \ge \gamma
l_1{ \prime d}  |Q|$.
\end{lemma}


{\bf \noindent Proof of Lemma \ref{t1-7}.} Consider the sets
$Q_0=Q$, $Q_{i}= 2Q_{i-1}$, for $i=1, \dots, s_1-h_1=s_2$, where
$h_1$ is the largest integer satisfying $2^{dh_1+d}< c$. If $Q_i$
was proper for all $i$, then $|Q_i|> \Volume (Q_i)$ and $\Volume
(Q_i) =2^d \Volume (Q_{i-1})$ and this  would imply that

\begin{equation}  |Q_{s_2}| > \Volume (Q_{s_2}) = 2^{ds_2} \Volume (Q_0) \ge
\frac{l_1^d}{2^{dk_1}} \frac{|Q|}{2^d} \ge  \frac{l^d
|Q|}{2^{dk_1+d}} \ge \frac{cm}  {2^{dk_1+d}} > m, \end{equation}

\noindent which is impossible as we assume $l_1Q \subset [m]$. (In
the second inequality we used the fact that $\Volume (Q_0) =
\Volume (Q) \ge \frac{|Q|}{2^d}$.)  Therefore, there is some $i$
between $1$ and $s_2$ for which $Q_i$ is not proper. Let $j$ be
the smallest such $i$. Thus, $Q_{j-1}$ is proper and $Q_j=2
Q_{j-1}$ is not. By Lemma \ref{t1-6}, there are constants $h_2$
and $\gamma_1$ such that $h_2 Q_{j-1}$ contains a proper GAP  $Q'$
of rank $(d-1)$ with cardinality at least $\gamma_1 |Q_{j-1}|$.
Without loss of generality we can assume that $h_2$ is a power of
2, $h_2=2^{h_3}$. By increasing $c$, we can assume that $h_1> h_3$
which guarantees  that $l_1'= h_2 2^{j} \le l_1$. The set $l_1' Q
= h_2 Q_{j-1}$ contains a proper GAP $Q'$ of rank $(d-1)$ and
cardinality

\begin{equation}  |Q'| \ge \gamma_1 |Q_{j-1}| = \gamma_1 2^{(j-1)d} |Q| \ge
\frac{\gamma_1}{h_2^d} l_1^{\prime d}  |Q| = \gamma l_1^{\prime d}
|Q|,
\end{equation}

\noindent where $\gamma =\frac{\gamma_1}{h_2^d}$, concluding the
proof. \hs

\subsection { Proof of Theorem \ref{1}}  Before
starting the proof, let us mention that all constants ($\gamma_1,
\gamma_2$ etc) in the proof depend on $d$, but do not depend on
$C$. By setting $C$ sufficiently large, we can satisfy all
relations required between these constants. Without loss of
generality, we can assume $l$ is a power of two, $l=2^s$, where
$s$ is sufficiently large.  Consider the set sequence $A_0=A$,
$A_{i+1} = 2A_i$. We first need the following fact, which asserts
that for some $i$ significantly smaller than $s =\log_2 l$, the
ratio $|A_i|/ |A_{i-1}|$ is not too large.

\begin{fact} \label{t1-8} There is some $i \le \frac{d+1}{d+3/2}s $ such that $|A_i| \le
2^{d+3/2} |A_{i-1}|$. \end{fact}

{\bf \noindent Proof of Fact \ref{t1-8}.} Assume otherwise, then

\begin{equation}  \label{ln} |A_{ \frac{d+1}{d+3/2}s}|  \ge 2^{(d+3/2)
\frac{d+1}{d+3/2}s} |A_0| = 2^{(d+1)s}|A| = l^{d+1}|A| \ge Cln,
\end{equation}

\noindent a contradiction as $A_{ \frac{d+1}{d+3/2}s}$ is a subset
of $[ln]$ ($C$ is set to be larger than 1). The proof of the claim
is completed. \hs

\vskip2mm

Let $s_1$ be the first index where $|A_{s_1+1 }|  \le 2^{d+3/2}
|A_{s_1}|$. Lemmas \ref{t1-2}, \ref{bilu} and \ref{full1} imply
that there are constants $g_1$ and $\gamma_1$ depending only on
$d$ such that $2^{g_1} A_{s_1}$ contains a proper GAP $Q$ of rank
$d+1$ and cardinality at least $\gamma_1 |A_{s_1}|$. By the
definition of $s_1$

$$|A_{s_1}| \ge 2^{(d+3/2)s_1} |A|, $$

\noindent so

$$ |Q| \ge \gamma_1 2^{(d+3/2)s_1} |A|. $$

\noindent By setting $C$  sufficiently large, we can assume that
$s$ is sufficiently large so that $s \ge s_1 + g_1$ (notice that
$s_1 \le \frac{d+1}{d+3/2} s$). This implies that
$\frac{l}{2^{{s_1} +g_1} } Q$ is a subset of $lA$. Next we apply
Lemma \ref{t1-7} to $Q$ with $m= ln$, $l_1 =
\frac{l}{2^{s_1+g_1}}$, and $d+1$ instead of $d$. In order to
verify the conditions of this lemma, observe that

\begin{equation}  l_1^{d+1} |Q| \ge \frac{l^{d+1}} {2^{(s_1+g_1)(d+1)}}
\gamma_1 2^{(d+3/2)s_1} |A| \ge \gamma_1 lcn  2^{s_1/2 -g_1(d+1)}.
\end{equation}

\noindent Again by assuming that  $C$ is large, we could guarantee
that the condition of Lemma \ref{t1-7} is met. Lemma \ref{t1-7}
implies that we have a proper GAP $Q' \subset l_1'Q =2^{s_1'} Q$
of rank $d$ with cardinality at least

\begin{equation}   \label{Q'} \gamma_2 2^{s_1'(d+1)} |Q| \ge \gamma_2 \gamma_1
2^{(d+1)(s_1+s_1')} 2^{s_1/2} |A| , \end{equation}

\noindent where $s_1+s_1' \le s$. The GAP $P = 2^{s-s_1-s_1'} Q'$
is a subset of $2^sA =lA$ and its volume is

\barray 2^{d( s-s_1-s_1') } \Volume (Q') \ge  2^{d( s-s_1-s_1')}
\frac{|Q'|}{2^d} &\ge& 2^{d( s-s_1-s_1'-1)}  \gamma_1 \gamma_2
2^{(d+1)(s_1+s_1')} 2^{s_1/2} |A| \\ &=& \gamma_1 \gamma_2 2^{ds}
|A| 2^{3s_1/2+s_1' -d} \\ &\ge& \frac{\gamma_1 \gamma_2}{2^d} l^d
|A|. \earray

\noindent Since $P$ has rank $d$, its longest edge forms an AP of
length at least $$\Big( \frac{\gamma_1 \gamma_2}{2^d} l^d |A|
\Big)^{1/d} =\Omega (l |A|^{1/d}), $$  completing the proof of
Theorem \ref{1}. \hs

\vskip2mm

{\bf \noindent Remark.}  The reader may notice that in this proof
we used the estimate on the cardinality of $Q'$, but did not use
the fact that $Q'$ is proper. The properness of $Q'$, however, is
critical in the next proof.

\subsection { Proof of Theorem \ref{111}} Without loss of generality, we can assume that
$0 \in A$. Consider $Q'$ as in the proof of Theorem \ref{1}. Again
by increasing $C$, we may assume that $s-s_1 -s_1'$ is lower
bounded by a sufficiently large constant. Consider the GAP $Q^{''}
= 2^{s-s_1-s_1'-g_2} Q'$, where $g_2$ is a large constant
satisfying $ s-s_1-s_1'-g_2 \ge 0$. Since $0 \in A$, $Q^{''}$ is a
subset of $lA$. Moreover, as $Q'$ and $Q^{''}$ are of ranked $d$,
we have, using inequality (\ref{Q'}), that

\barray \Volume (Q^{''}) &\ge&  2^{(s-s_1-s_1'-g_2)d}   \gamma_2
\gamma_1 2^{(d+1)(s_1+s_1')} 2^{s_1/2} |A| \\ &=& \gamma_1
\gamma_2 2^{sd} 2^{\frac{3s_1}{2} + s_1' -g_2 d} |A| \\ &=& \Omega
(l^d |A|). \earray

\noindent We are going to examine  two cases:

\vskip2mm

{\bf \noindent Case 1:} $Q^{''}$ is proper. In this case $lA$
contains the  proper GAP $Q^{''}$  of rank $d$ and volume $\Omega
(l^d |A|)$. So we are done by setting $d'=d$.

\vskip2mm

{\bf \noindent Case 2:} $Q^{''}$ is not proper. Now we make a
crucial use of the fact that $Q'$ is proper. The properness of
$Q'$ implies that  there is a positive integer  $s_2 \le
s-s_1-s_1'-g_2 \le s$ such that $\frac{1}{2^{s_2}} Q^{''}$ is
proper. As usual, we choose $s_2$  to be the smallest such an
integer, which implies that  $\frac{1}{2^{s_2-1}} Q^{''}=
\frac{2}{2^{s_2}} Q^{''}$ is not proper. Applying Lemma \ref{t1-6}
to $\frac{1}{2^{s_2}} Q^{''}$ we obtain a GAP $Q^{'''}$ of rank
$d-1$ and volume

$$\Omega (\Volume (\frac{1}{2^{s_2}} Q^{''})) =\Omega
(\frac{1}{2^{ds_2}} \Volume (Q^{''}) =\Omega (\frac{1}{2^{ds_2}}
l^d |A|). $$

\noindent Furthermore, there is a constant $g_3$ such that
$Q^{''''} = 2^{s_2-g_3} Q^{'''}$ is a subset of $lA$. The GAP
$Q^{''''}$ has rank $d-1$ and volume

$$2^{(s_2-g_3)(d-1)} \Volume (Q^{'''}) =\Omega (2^{s_2(d-1)}
\frac{1}{2^{ds_2}}l^d |A|) = \Omega (2^{-s_2} l^d |A|). $$

\noindent Since $s \ge s_2$, $2^{-s_2} \ge 2^{-s}= l^{-1}$. Thus,
the volume of $Q^{''''}$ is  $\Omega (l^{d-1} |A|)$. Now if
$Q^{''''}$ is proper then we are done by setting $d'=d-1$.
Otherwise we repeat the analysis of Case 2 to obtain a GAP of rank
$d-2$ and so on. This repetition cannot continue forever so sooner
or later we must obtain a proper GAP of some rank $d' <d$ which
satisfies the claim of the theorem. \hs

\section {Sums of different sets}

The goal of this section is to  generalize the results in Section
3 by considering the sum of different sets, instead of the sum of
the same sets. Given $l$ sets $A_1, \dots, A_l$, we define

$$A_1+ \dots +A_l =\{a_1 +\dots + a_l| a_i \in A_i,  1\le i \le l
\}. $$

\noindent We obtain the following generalization of Theorem
\ref{111}.

\begin{theorem}  \label{211} For any fixed positive integer $d$
there are positive constants $C$ and $c$ depending on $d$
 such that the following holds. Let $A_1, \dots, A_l$ be subsets  of $[n]$ of size $|A|$ where
 $l$ and $|A|$ satisfy $l^d |A| \ge Cn$.
Then $A_1 +\dots +A_l$ contains a GAP of rank $d'$ and volume at
least $cl^{d'} |A|$, for some integer $1 \le d' \le d$.
\end{theorem}

\noindent The following corollary generalizes Theorem \ref{1}.

\begin{corollary} \label{212} For any fixed positive integer $d$
there are positive constants $C$ and $c$ depending on $d$
 such that the following holds. Let $A_1, \dots, A_l$ be subsets  of $[n]$ of size $|A|$ where
 $l$ and $|A|$ satisfy $l^d |A| \ge Cn$.
Then $A_1 +\dots +   A_l$ contains an arithmetic progression of
length $cl|A|^{1/d}$.
\end{corollary}

Corollary \ref{211} has a  nice application. In Section 6, we use
this corollary  to confirm a conjecture of Folkman posed in 1966.

\subsection {The basic idea}

The basic idea behind the proof of Theorem \ref{211} is the
following. Given the sets $A_1, \dots A_l$ as in Theorem
\ref{211}, we are going to show that there are numbers $l', n'$
and a set $A'$ such that

\begin{itemize}

\item $l'A'$ is a subset of $ A_1 +\dots +A_l$; $A'$ is a subset
of $[n']$.

\item   $l', n', |A'|$ satisfy the conditions of Theorem
\ref{111}. (This can be done by setting the constant $C$ in
Theorem \ref{211} much larger than the constant $C$ in Theorem
\ref{111}.)

\item $(l')^{d'} |A'| =\Omega ( l^{d'} |A|)$ for all $1 \le d' \le
d$.

\end{itemize}

\noindent An application of Theorem \ref{111} to the triple $(l',
n', A')$ immediately  implies the statement of Theorem \ref{211}.

\vskip2mm

The proof of Theorem \ref{211} uses a technical lemma, Lemma
\ref{Y} below. This lemma provides a sufficient condition for the
existence of a sumset of form $l'A'$ in a sumset of different
sets. The verification of this lemma requires extensions of the
filling lemmas described in Section 4. These extensions is the
topic of the next subsection.

\subsection {Filling with different sets}

In the proof of Theorem \ref{211}, we shall need the following
lemma, which generalizes Lemma \ref{full10} the way Theorem
\ref{211} generalizes Theorem \ref{1}. This lemma was proved in an
earlier paper. For the readers' convenience, we include the proof
in Appendix A.

\begin{lemma} \label{sum} For any positive constant $\gamma$ and positive integer $d$,
there is a positive constant $\gamma'$ and a positive integer $g$
such that the following holds. If $X_1, \dots, X_g$ are subsets of
a generalized  arithmetic progression $P$ of rank $d$ and $|X_i|
\ge \gamma { \hbox{Vol} }(P)$ then $X_1 +\dots + X_g$ contains a
generalized  arithmetic progression $Q$ of rank $d$ and
cardinality at least $\gamma' {\hbox{Vol}} (P)$. Moreover, the
distances of $Q$ are multiplies of the distances of $P$.
\end{lemma}

\noindent One can further strengthen this lemma by requiring $Q$
be proper. The proof is similar to the proof of proper filling
lemma, Lemma \ref{full1}.

\begin{lemma} \label{sum1} For any positive constant $\gamma$ and positive integer $d$,
there is a positive constant $\gamma'$ and a positive integer $g$
such that the following holds. If $X_1, \dots, X_g$ are subsets of
a generalized  arithmetic progression $P$ of rank $d$ and $|X_i|
\ge \gamma { \hbox{Vol} }(P)$ then $X_1 +\dots + X_g$ contains a
proper generalized  arithmetic progression $Q$ of rank $d$ and
cardinality at least $\gamma' {\hbox{Vol}} (P)$. Moreover, the
distances of $Q$ are multiplies of the distances of $P$.
\end{lemma}

Later on, we shall refer to Lemmas \ref{sum} and \ref{sum1} as the
general filling and general proper filling lemmas, respectively.

\subsection { The main lemma of Theorem \ref{211}}

We are now in position to present and prove the main lemma of the
proof of Theorem \ref{211}.

\begin{lemma} \label{Y} For every positive constant $c$ there are
positive constants $\ep$ and $d$ depending on $c$ such that the
following holds. If the sets $X_1, \dots, X_l$, each of
cardinality $|X|$, satisfy $|X_1 +X_i | \le c|X|$ for all $2 \le i
\le l$, then there is a proper GAP $Q$ of rank at most $d$ and
cardinality at least $\ep |X|$ and a number $l' \ge \ep l$ such
that the sum $X_1+\dots +X_l$ contains a translation of $l'Q$.

\end{lemma}

{\bf \noindent Proof of Lemma \ref{Y}.} The condition $|X_1 +X_i|
\le c|X|$ and  Freiman's theorem imply that  $X_1$ is contained in
a GAP $R$ with constant rank and volume $O(|X|)$. Consider $X_i$,
for some $2 \le i \le l$. We say that two elements $x$ and $y$ of
$X_i$ are equivalent if $x-y \in R-R$. It is trivial that if $x$
and $y$ are not equivalent then $x+X_1$ and $y+X_1$ are disjoint
sets. Since $|X_1 +X_i| \le c|X|$ where $|X|=|X_1|=|X_i|$, the
number of equivalent classes is at most $c$. It follows that there
is a class with cardinality $\Omega (|X|)$; let us call this class
$Y_i$. $Y_i$ is a translation of a subset $Z_i$ (of constant
density)  of $R$. The hidden constants in the asymptotic notations
depend on $c$.

\vskip2mm

Consider the sets $Z_2, \dots, Z_l$. These sets are subsets of $R$
and $|Z_i| \ge \gamma \Volume (R)$ for some positive constant
$\gamma$ depending on $c$. Let $g$ be  a large constant integer.
With the exception of at most $g-1$ sets,  we partition the
$Z_i$'s into $l_1 =\lfloor (l-1)/g \rfloor$ disjoint groups of
size $g$: $G_1, \dots, G_{l_1}$. Thus each group $G_j$ contains
$g$ sets, each of which is a subset of $R$ with cardinality
$\gamma \Volume (R)$ for some positive constant $\gamma$. By
setting $g$ sufficiently large, the general  filling lemma (Lemma
\ref{sum}) applies and shows that the sum of the sets in any group
$G_j$ contains a proper GAP $Q_j$ of cardinality $\Omega (\Volume
(R))$. Moreover, the rank of $Q_j$ is the same as the rank of $R$
and the differences of $Q_j$ are multiples of the differences of
$R$.

\vskip2mm Since $|Q_j| =\Omega (\Volume (R))$, there are only
$O(1)$ choices for the difference set of  $Q_j$ (for the
definition of difference sets, see Section 2). Thus, a constant
fraction of the $Q_j$'s has the same difference set. Without loss
of generality, we may assume that these $Q_j$'s are $Q_1, Q_2,
\dots, Q_{l_2}$, where $l_2 =\Omega (l_1)$. \vskip2mm

Since $|Q_j| =\Omega (\Volume (R))$, the length of the $h$th edge
of $Q_j$ is $\Omega(1)$ times the length of the corresponding edge
of $R$, for all $1 \le h \le \rank (R)$. Thus the lengths of the
$h$th edge of the $Q_j$'s are within a constant factor from each
other, for all $1\le k \le \rank (R)$. This implies that the
intersection of the boxes $B_{Q_1}, \dots, B_{Q_{l_2}}$ contains a
box $B$ with volume $\Omega (\Volume (R))$ (for the definition of
these boxes, see subsection 3.13). Let $m_1, \dots, m_d$ be the
lengths of the edges of $B$ and $(a_1, \dots, a_d)$ be the
(common) set of differences of $Q_1, \dots, Q_{l_2}$. It follows
that each of $Q_1, \dots, Q_{l_2}$ contains a translation of the
proper GAP $Q= \{ a_1 x_1 +\dots a_d x_d| 0 \le x_i \le m_i \}$
($Q$ is proper because the $Q_j$'s are so). We have that

$$|Q| =|B| = \Omega (\Volume (R)) = \Omega (|X|), $$

\noindent and

$$l_2 =\Omega (l_1) = \Omega (l).$$

\noindent Moreover, a translation of $l_2 Q$ is contained in  $Q_1
+\dots Q_{l_2}$ and a translation of $Q_1 +\dots Q_{l_2}$ is
contained in $X_1 +\dots + X_l$. So $X_1 +\dots +X_l $ contains a
translation of $l_2Q$, completing the proof. \hs

\subsection {Proof of Theorem \ref{211}}
With the main lemma in hand, we are ready to conclude the proof of
Theorem \ref{211}. In order to find a triplet $(A', l', n')$ as
desired, we are going to apply the so-called tree argument. This
argument was introduced in \cite{szemvu2} and, in spirit, works as
follows. Assuming that we want to add several sets $A_1, \dots,
A_l$. We shall add them in a special way following an algorithm
which assigns sets to the vertices of a tree. A set of any vertex
contains the sum of the sets of its children. If the set at the
root of the tree is not too large, then there is a level where the
sizes of the sets do not increase (compared with the sizes of
their children)  too much. Thus, we can apply Freiman's inverse
theorems at this level to deduce useful information. The creative
part of this argument is to come up with a proper algorithm which
suits our need.

\vskip2mm

The reader has already met a simple version of this argument in
the proof of Theorem \ref{1}. In that proof, the sets at the
leaves of the tree are copies of $A$, the sets at a level $i$ are
copies of $2^iA$ and the set at the root is $lA$. A set of any
vertex is the sum of the sets at its two children.

\vskip2mm

\vskip2mm  The algorithm in the current case is more complicated.
Before describing it,  let us assume, without loss of generality,
that $l$ is a power of 4 ($l=4^s$) and $|A_i|=n_1$ and $0 \in A_i$
for all $1 \le i \le l$. Set $A_i= A_1^1$ for $i=1, \dots, l$ and
$l_1=l$. Here is the description of the algorithm.

\vskip2mm {\noindent \it The algorithm.} At the $t^{th}$ step, the
input is a sequence $A_1^{t}, \dots, A_{l_t}^{t}$ of the same
cardinality $n_t$ where $l_t$ is an even number. Choose a pair $1
\le i < j \le l_t$ which maximize $|A_i^{t} + A_j ^{t}|$ (if there
are many such pairs choose an arbitrary one).  Denote the sum
$A_i^t + A_j^t$ by $A'_1$.  Remove $i$ and $j$ from the index set
and repeat the operation to obtain $A_2'$ and so on. After $l_t/2$
operations we obtain  a sequence $A_1', \dots, A_{l_t/2}'$ of sets
with decreasing cardinalities. Define $l_{t+1} = l_t/4$. Consider
the sequence $A'_1, \dots, A'_{l_{t+1}}$ and truncate all but the
last set  so that all of them have the same cardinality (which is
$|A'_{l_{t+1}}|$). The truncated sets will be named $A_1^{t+1},
\dots, A_{l_{t+1}}^{t+1}$ and they form the input of the next
step.  It is clear  that $l_t =\frac{l_1}{4^{t-1}} $ for all
plausible $t$'s. The algorithm halts at time $s+1$ where $l_{s+1}
=1$.

\vskip2mm

Notice that $A^{t+1}_{l_{t+1}}$ is a subset of $[2^t n]$, so
$n_{t+1} \le 2^tn$. We first show that there is some $t\le s$ so
that $n_{t+1} \le 4^{d+1} n_t$. Assume otherwise. Then

$$n_{s+1} \ge (4^{d+1})^s n_1 = (4^s)^d 4^s |A| = 4^s l^d |A| >
2^s n, $$

\noindent a contradiction. In the following, let $t$ be the first
index so that $n_{t+1} \le 4^{d+1} n_t$. By the description of the
algorithm, there are $l_t/2$ sets among the sets $A^t_i$'s such
 that every pair of them have cardinality at most $n_{t+1} \le
4^{d+1} n_t$. Let us call these sets $B_1, \dots, B_{l_t/2}$. We
have

\begin{itemize}

\item $ |B_1| =\dots =B_{l_t/2} = n_t,$

\item $B_i$'s are subsets of the interval $[2^{t-1} n]$,

\item $|B_i +B_j| \le 4^{d+1} n_t$, for all $1\le i < j \le
l_t/2$.

\end{itemize}

\noindent By Lemma \ref{Y}, the sum $B_1 + \dots + B_{l_t/2}$
contains a translation of $l'A'$, where $l' \ge \ep l_t/2$ and
$A'$ is a proper GAP with cardinality at least  $\ep n_t$ and
$\ep$ is a positive constant depending on $d$.  Moreover, $A'$ is
a subset of $[k_1 2^{t-1} n]$, for some constant $k_1$ depending
on $d$. Set $n' =k_12^{t-1} n$. To conclude the proof, let us
verify that $l',n'$ and $A'$ satisfy the required relations. First
of all

\begin{equation}  \label{l'} (l')^d |A'| \ge (\ep l_t/2)^d \ep n_t \ge
\frac{\ep^{d+1}}{2^d} \frac{l^d}{4^{(t-1)d}} 4^{(d+1) (t-1)} |A|
\ge \frac{\ep^{d+1}}{2^d} 2^{t-1} l^d |A|. \end{equation}

\noindent Since $l^d |A| \ge Cn$, it follows that

$$(l')^d |A'| \ge \frac{C \ep^{d+1}}{2^d}  2^{t-1}n = \frac{C \ep^{d+1}}{k_1 2^d}
n'. $$

\noindent By increasing $C$ (notice that $\ep$ and $k_1$ do not
depend on $C$), we can assume that $(l')^d|A'|/n' $ is
sufficiently large. This guarantees that the condition of Theorem
\ref{111} is met. Replacing $d$ by $d'$ in  (\ref{l'}) one can
verify that for any $d' \le d$

$$(l')^{d'} |A'| =\Omega (l^{d'} |A|),$$

\noindent  concluding the proof. \hs

\section { Folkman's conjecture on subcomplete sequences}

For a  (finite or infinite) set $A$,  $S_A$ denotes the collection
of subset sums of $A$

$$S_A=\{ \sum_{x\in B} x| B \subset A, |B| < \infty \}. $$

An infinite sequence $A$ of positive integers is {\it subcomplete}
if $S_A$ contains an infinite arithmetic progression. Subcomplete
sequences have been studied extensively and we refer the reader to
Section 6 of the monograph \cite{EG} by Erd\"os and Graham for a
survey. For an infinite sequence $A$, we use $A(n)$ to denote the
number of elements of $A$ between 1 and $n$. This number could be
larger than $n$ as $A$ might contain the same number many times.
In 1966, Folkman made the following conjecture

\begin{conj}  \label{folk} There is a constant $C$ such that the following holds. If $A = \{a_1 \le a_2 \le a_3 \le \dots \}$ is an
infinite non-decreasing sequence of  positive integers and $A(n)
\ge Cn$, for all sufficiently large $n$, then $A$ is subcomplete.
\end{conj}

Folkman's conjecture was considered by Erd\"os and Graham as the
most important conjecture concerning subcomplete sequences
(\cite{EG}, Section 6). Folkman himself proved that the conjecture
holds under a stronger condition that $A(n) \ge n^{1+\ep}$, where
$\ep $ is an arbitrary positive constant. The conjecture is sharp,
as one cannot replace $n$ by $n^{1-\ep}$. To show this, let us
present an observation of Erd\"os \cite{erd}.

\begin{fact} \label{erdos} Consider an infinite sequence $A =
\{a_1, a_2, \dots\}$. If

\begin{equation}  \label{limsup} \limsup_{i \rightarrow \infty} ( a_i
-\sum_{j=1}^{i-1} a_j) \rightarrow \infty, \end{equation}

\noindent then $A$ is not subcomplete. \end{fact}

 To
verify Fact \ref{erdos}, notice that if $A$ is subcomplete and $d$
is the difference of an infinite arithmetic progression contained
in $S_A$, then $d$ is lower bounded by $\limsup_{i \rightarrow
\infty}  (a_i -\sum_{j=1}^{i-1} a_j)$.

\vskip2mm

For any fixed $\ep >0$, it is simple to find a non-decreasing
sequence $A$ such that $A(n) =\Omega (n^{1-\ep})$ and $A$
satisfies (\ref{limsup}).

\vskip2mm  Using a special case of Theorem 5.1 (Corollary
\ref{212}), we are able to confirm Folkman's conjecture.

\begin{theorem}  \label{folkman} There is a constant $C$ such that the following holds.
If $A = \{a_1 \le a_2 \le a_3 \le \dots \}$ is an infinite
non-decreasing sequence of  positive integers and $A(n) \ge Cn$,
for all sufficiently large $n$, then $A$ is subcomplete.
\end{theorem}

The rest of this section is devoted to the proof of Theorem
\ref{folkman}, which relies on Corollary \ref{212}. First, we
prove a sufficient condition for subcompleteness. This condition
is of independent interest and will be used for another problem in
Section 9.  To complete the proof, we show that any sufficiently
dense sequence satisfies this condition. This part of the proof
makes a significant use of Corollary \ref{212}.

\subsection {A sufficient condition for subcompleteness}

We say that a sequence $A$ admits a {\it good} partition if it can
be partitioned into two subsequences $A'$ and $A^{''}$ with the
following two properties

\begin{itemize}

\item There is a number $d$ such that   $S_{A'}$ contains an
arbitrary long arithmetic progression with difference $d$.

\item Let $A^{''} = b_1 \le b_2 \le b_3 \le \dots$. For any number
$K$, there is an index $i(K)$ such that $\sum_{j=1}^{i-1} b_j \ge
b_i +K$ for all $i \ge i(K)$.

\end{itemize}

\noindent Admitting a good partition is a  sufficient condition
for  subcompleteness.

\begin{lemma} \label{subcomplete} Any sequence $A$ which admits a
good partition is subcomplete. \end{lemma}

{\bf \noindent Proof of Lemma \ref{subcomplete}.} We start with a
definition.

\begin{defn} \label{net} An infinite sequence $B =\{b_1 \le b_2 \le b_
3 \le \dots \}$ is a $(d,L)$-net if $b_{i+1} -b_i < L$ and is
divisible by $d$ for all $i =1,2 \dots.$ \end{defn}

Observe that if $B$ is a $(d, L)$-net and $Q$ is a finite
arithmetic progression with difference $d$ and  length larger than
$L/d$, then $B+Q$ contains an infinite  arithmetic progression
with difference $d$. This observation is the leading idea in what
follows.

\vskip2mm

Assume that $A$ admits a good partition and let $Q_0, Q_1,  Q_2,
\dots$ be  arithmetic progressions
 with the same difference $d$ and
strictly increasing lengths contained in $S_{A'}$. The existence
of the $Q_i$'s is guaranteed by the first property of a good
partition. \vskip2mm

Next, we focus on $A^{''}$. Let $X$ be the set of divisors $d'$ of
$d$ with the following property: All but at most finite elements
of $A^{''}$ are divisible by $d'$. Since $1 \in X$, $X$ is not
empty and thus has a maximum element $d_1$. By throwing away
finite elements, we can assume that all  elements of $A^{''}$ are
divisible by $d_1$. Next, discard all elements $y$ (in the
remaining sequence) if there is only a finite number elements of
$A^{''}$ which equal $y$ modulo $d$. Again, we discard only a
finite number of elements so the remaining sequence still has the
same density as $A^{''}$. Thus, we can assume that $A^{''} = \{b_1
d_1 \le b_2 d_1 \le \dots \}$ where the $b_i$'s have the following
property:
 Let
$b_i'$ be the remainder when dividing $b_i$ by $d$.  For each $i$,
there are infinitely many $j$'s such that $b_i' = b_j'$. Moreover,
the greatest common divisor of the $b_i'$'s equals one modulo $d$
by the definition of $d_1$. We next need the following elementary
fact, which is a consequence of the Chinese remainder theorem.

\begin{fact} \label{mod} Let  $1 \le z_1 \le z_2 \le \dots \le z_h < d$
be positive integers. If $\gcd (z_1, \dots, z_h)=1 (\hbox{mod}
\,\, d)$, then there are integers $0 \le a_1 , \dots, a_h < d$
such that $\sum_{j=1}^h a_j z_j \equiv 1 (\hbox{mod} \,\, d)$.
\end{fact}

By Fact \ref{mod} and the property of $A^{''}$ described in the
previous pararaph, we can find $(d-1)$ mutually disjoint finite
subsets $X_1, \dots, X_{d-1}$ of $A^{''}$ so that the sum of the
elements in each subset equals $d_1$ modulo $d$. Denote these sums
by $x_1 d+ d_1, \dots, x_{d-1} d+ d_1$, where $x_i$'s are
non-negative integers. For any arithmetic progression $Q_j$ with
length $l \ge 3 (x_1 + \dots + x_{d-1})$, the set $Q_j + S_{
\{x_1d+d_1, \dots, x_{d-1} d+ d_1 \}}$ contains an arithmetic
progression with difference $d_1$ and length at least $l/2$
(recall that $Q_j$ has difference $d$ which is divisible by
$d_1$). Since the lengths of the $Q_j$'s go to infinity with $j$,
we can conclude that $S_{A'} + S_{ \{x_1d+d_1, \dots, x_{d-1} d+
d_1\}}$ contains an arbitrarily long arithmetic progression with
difference $d_1$.

\vskip2mm Set  $A^{'''}= A^{''} \backslash \cup_{i=1}^{d-1} X_i$;
to complete the proof of the lemma, it suffices to prove that
$S_{A^{'''}}$ contains a $(d_1, L)$-net for some constant $L$. Let
$S_{A^{'''}} =\{s_1 < s_2 < \dots \}$. Every elements of $A^{''}$
is divisible by $d_1$ and so are all $s_i$'s. Therefore, it
suffices to exhibit the existence of a constant $L$ satisfying
$s_{i+1} -s_i \le L$ for all $i$. The existence of $L$ follows
directly from the following observation, due to Graham \cite{Gra},
and the second property of a good partition (this is the only
place where we use this property).


\begin {fact} \label{graham} Let $Y= y_1 < y_2 < \dots $ be an infinite
sequence of positive integers and $S_{Y} = \{s_1 < s_2 < \dots
\}$. If $y_{m+1 } \le \sum_{i=1}^m y_i$ for all sufficiently large
$m$, then there is some $L$ such that $s_{i+1} -s_i \le L$ for all
$i$. \end{fact}

\noindent Fact \ref{graham} is not too hard and the reader might
want to consider it as an exercise.

\subsection{Proof of Theorem \ref{folkman}}

\noindent We first present a lemma  which provides a link between
good partitions and subcompleteness. This lemma is a simple, but a
bit tricky,  consequence of  Corollary \ref{212}.

\begin{lemma} \label{folkman1} There is a constant $C$ such that the following holds. If $A$ is a multi-set of positive
integers between 1 and $n$ and $|A| \ge Cn$, then $S_A$ contains
an arithmetic progression of length $n$. \end{lemma}

{\bf \noindent Proof of Lemma \ref{folkman1}.} We show that the
same constant $C$ in Corollary \ref{212} suffices. Without loss of
generality, we assume that $C$ is an integer and $|A|=Cn$. If the
multi-set $A$ contains an element $a$ of multiplicity $n$, then
the arithmetic progression $a, 2a, \dots , na$ is a subset of
$S_A$ and we are done. In the other case, we can partition the
$Cn$ elements of $A$ into $n$ sets $X_1, \dots, X_n$ such that
each $X_i$ consists of  exactly $C$ different elements. The sum
$X_1 +\dots +X_n$ is a subset of $S_A$. Corollary \ref{212}
implies that the sum $X_1 +\dots +X_n$ contains an arithmetic
progression of length $n$, given that $C$ is sufficiently large.
This concludes the proof of the lemma. \hs

\vskip2mm

With Lemma \ref{folkman1} in hand, we are in a position to prove
that the sequence $A$ in Theorem \ref{folkman} admits a good
partition, provided that the constant $C$ in this theorem is
sufficiently large. The partition is the most natural one. Assume
that the elements of $A$ are ordered non-decreasingly $A=a_1 \le
a_2 \le a_3 \le \dots$; $A'$ consists of the elements with odd
indices, $A^{''}$ consists of those with even indices.

\vskip2mm

  By
definition, $A^{''} = \{a_2, a_4, a_6, \dots\}$. Since  $A(n)\ge
Cn$ for all sufficiently large $n$ (recall that $A(n)$ is the
number of elements of $A$ between 1 and $n$), for every
sufficiently large even number $j$,

$$a_j \le j/C \le j/5 \le a_2 + a_4 + \dots + a_{j-2} - j/4, $$

\noindent which guarantees the property required for $A^{''}$.

\vskip2mm It remains to check the property  concerning $A'$. As
$A$ has density $Cm$, $A'$ has density $Cm/2$ so we can assume
that $A' = \{ b_1 \le b_2 \le \dots \}$, where $b_m \le 2m/C$ for
all sufficiently large $m$. Let $A'[m] $ be the set consisting of
the first $m$ elements of $A'$. Fix a sufficiently large $m$ and
define $A_0 =A'[m]$ and $A_i = A'[2^i m] \backslash A '[2^{i-1}
m]$. The set $A_i$ has $2^{i-1} m$ elements and is a subset of the
interval $[2^{i+1}m/C]$.

\vskip2mm

To conclude the proof, we make use of the following lemma, proved
in \cite{szemvu}

\begin{lemma} \label{sumtwo1}  Let $P$ be a generalized arithmetic progression of rank
2, $P =\{x_1 a_1 + x_2a_2 | 0\le x_i \le l_i \}$, where  $l_i \ge
5a_{3-i}$ for $i=1,2$. Then $P$ contains an arithmetic progression
of length $l_1+l_2$ whose difference is $\gcd(a_1,a_2)$.
\end{lemma}

\vskip2mm

By Lemma  \ref{folkman1} (provided that $C$ is sufficiently
large), $S_{A_i}$ contains an  arithmetic progression $P_i$ of
length $l_i = 2^{i+1}m/C$ for all $i$. Set $Q_0=P_0$ (and assume
that $d_0$ is the difference of $Q_0$) and consider the
generalized arithmetic progression $Q_0 +P_1$. This is a
generalized arithmetic progression of rank 2 with volume $l_0
l_1$. Moreover, this two dimensional generalized arithmetic
progression  is a subset of an interval of small lenght, so one
can easily check that its differences are relatively small and
satisfy the assumption of Lemma \ref{sumtwo1}. This lemma implies
that $Q_0+P_1= P_0 + P_1$ contains an arithmetic progression $Q_1$
of length $l_0 + l_1-2 $ with difference $d_1$ which is a divisor
of $d_0$. (The $-2$ term comes from the fact that in Lemma
\ref{sumtwo1}, the edges of $P$ have length $l_1+1$ and $l_2+1$,
respectively; this term, of course, plays no role.) Similarly, by
considering $Q_1 +P_2$ we obtain an arithmetic progression $Q_2$
of length $l_0+l_1 + l_2-3$ with difference $d_2$ which is a
divisor of $d_1$ and so on. The sequence $d_0, d_1, d_2, \dots$ is
non-increasing, so there is an index $j$ so that $d_i =d_j=d$ for
all $i \ge j$. The arithmetic progressions $Q_j, Q_{j+1}, Q_{j+2},
\dots$ have strictly increasing lengths and the same difference
$d$. Moreover, each $Q_i$ is a subset of $S_{A'}$ and this
completes the proof. \hs

\section {Sumsets with   distinct summands}

In this section, we strengthen Theorem \ref{111} in another
direction. Instead of the sumset $lA$, we are going to consider
the much more restricted sumset $l^{\ast} A$, which consists of
the sums $a_1 +\dots +a_l$ where the $a_i$'s are different
elements of $A$.

\begin{theorem}  \label{311} For any fixed positive integer $d$
there are positive constants $C$ and $c$ depending on $d$
 such that the following holds. For any positive integers $n $
and $l$ and any set  $A \subset [n]$ satisfying $l \le |A|/2$ and
$l^{d} |A| \ge C n$,  $l^{\ast} A$ contains a proper GAP of rank
$d'$ and volume at least $cl^{d'} |A|$, for some integer $1 \le d'
\le d$.
\end{theorem}

\vskip2mm The requirement that the summands must be different
usually poses a great challenge in additive problems. One of the
most well-known examples is the celebrated Erd\"os-Heilbronn's
conjecture. In order to describe this conjecture, let us start
with the classical Cauchy-Danveport theorem which asserts that if
$A$ is a set of residues modulo $n$, where $n$ is a prime, then
$|2A| \ge \min \{n, 2|A|-1 \}$. For $A$ being an arithmetic
progression, the bound is sharp. Now let us consider $2^{\ast} A$.
We want to bound $|2^{\ast} A|$ from below with something similar
to Cauchy-Danveport's bound. Observe that in the special case when
$A$ is an arithmetic progression, $|2^{\ast} A| = \min \{n,
2|A|-3\}$. Thus one may guess that

\begin{equation}  \label{EH} |2^{\ast} A| \ge \min \{n, 2|A|-3\} , \end{equation}

\noindent holds for any set $A$. This is what Erd\"os and
Heilbronn conjectured. While Cauchy-Davenport's theorem is quite
easy to prove, Erd\"os -Heilbronn's conjecture had been open for
about thirty years until it was solved by da Silva and Hamidoune
in 1994 \cite{SH}.

\vskip2mm It is now not so big a surprise that Theorem \ref{311}
is harder and deeper than both Theorem \ref{111} and Theorem
\ref{211}. The proof of Theorem \ref{311} uses Theorem \ref{111}
as a lemma and requires lots of additional  arguments, but let us
take a gentle start by introducing some simple ideas.

\subsection {The initial  ideas} \label{perfect}

The initial ideas in the proof of Theorem \ref{311} are similar to
those in the proof of Theorem \ref{211}.   We want to show that
there are numbers $l', n'$ and a set $A'$ such that

\begin{itemize}

\item  $A'$ is a subset of $[n']$ and  $l', n', |A'|$ satisfy the
conditions of Theorem \ref{111}, namely $l{\prime d} |A'|/n'$ is
sufficiently large.

\item $(l')^{d'} |A'| =\Omega ( l^{d'} |A|)$ for all $1 \le d' \le
d$.

\end{itemize}

\noindent In the rest of the proof, we call a triple $(A',l',n')$
{\it perfect} if it satisfies the above two conditions. If we
could  show that  there is a perfect triple $(A',l',n')$ such that
$l'A'$ is a subset of $ l^{\ast} A$, then an application of
Theorem \ref{111} to this triple immediately implies the statement
of Theorem \ref{311}.

\vskip2mm It is useful to notice that in Theorem \ref{311},
instead of the assumption $l \le |A|/2$, we can afford a stronger
assumption that $l \le \ep |A|$ for any positive constant $\ep$,
at the cost of increasing the constant $C$. One can argue as
follows. First one puts aside $(1-\ep)l$ elements from $A$. Next,
consider the pair $(A_1, l_1)$ where $A_1$ is the set of the
remaining $|A|-(1-\ep)l$ elements and $l_1 = \ep l$. It is trivial
that $l_1 \le \ep |A_1|$. On the other hand, the sum of an element
from $l_1^{\ast}A_1$ and the sum of the $(1-\ep)l$ elements put
aside is an element of $l^{\ast}A$. So if $l_1^{\ast}A_1$ contains
a proper GAP $P$, then $l^{\ast}A$ contains a translation of $P$.

\vskip2mm The above argument also shows that  for any $l_1  <l$,
if $A_1$ is a subset of at most  $|A|-(l-l_1)$ elements of $A$,
then $l_1 ^{\ast} A_1$ is a subset of a translation of $l^{\ast}
A$.

\vskip2mm In the proof of Theorem \ref{311}, we shall assume that
$l \le \ep|A|$, whenever needed. We shall also assume that $l-1$
elements of $A$ are put aside in case we need them to create the
sum of exactly $l$ elements. These assumptions provide  us some
flexibility in constructing a perfect triple. In particular, we
shall not need to show that $l'A'$ is a subset of $l^{\ast}A$; it
suffices to show that $l'A'$ is  a subset of a translation of
$\tilde l ^{\ast} A$, for some $\tilde l \le l$.

\vskip2mm The main part of the proof is to construct a perfect
triple and this is significantly harder than what we did in the
proof of Theorem \ref{211}. However, when $|A|$ is large the
construction is relatively simple and we start with this case. The
treatment of the harder case when $|A|$ is relatively small starts
in subsection 7.5, where we present a key structural lemma. The
proof of this lemma occupies the rest of this section. In the next
section, Section 7, we present the rest of the proof of  Theorem
\ref{311}.

\subsection { The case when $|A|$ is large}  \label {largeA}

Let $A_1$ be a subset of $A$ with cardinality $l-1$ and  set $A_2=
A\backslash A_1$. Since $|A| \ge 2l$,

\begin{equation}  \label{theorem20} |A_2| \ge \frac{|A|}{2}.\end{equation}

\vskip2mm
 We assume, with foresight (and with room to spare), that $|A|^2 \ge
 80C
n \log_2 n$ and $l^d|A| \ge 160 \times 2^d Cn$, where $C$ is the
constant in Theorem \ref{111}.

\vskip2mm

Define $m_i= 2^i$ for all $ 1\le i \le t$, where $t$ is the
smallest index such that $m_t \ge |A_2|/2$. Since $|A_2| \le |A|
\le n$,  $t \le \log_2 n$. Let $S_i$ be the set of those numbers
in $[2n]$ which can be represented as the sum of two different
elements in $A_2$ in at least $m_i$ and less than $m_{i+1}$ ways.
It is essential to observe that  $m_iS_i$ is a subset of
$(2m_i)^{\ast} A$. On the other hand, a simple double counting
argument gives

\begin{equation}  \label{theorem21} \sum_{i=1}^t m_i |S_i| \ge { {|A_2|}
\choose 2} -4n \ge q= \frac{|A|^2}{5}. \end{equation}

Next, we split $\sum_{i=1}^t m_i |S_i|$ into three parts. The
first part comprises those $m_i |S_i|$ where $m_i|S_i| \le
\frac{q}{4t}$. Obviously, the contribution of this part to the sum
is at most $t \frac{q}{4t} =\frac{q}{4}$. The second part consists
of those $m_i |S_i|$ where $|S_i | \le \frac{|A_2|} {40}$. Since
the sequence $m_i$ is geometric, the sum of all $m_i$'s is bounded
from above by $2|A_2|$. Thus, the contribution of the second part
is upper bounded by $2|A_2| \frac{|A_2|} {40}\le \frac{q}{4}$. The
third part contains the remaining $m_i|S_i|$'s and, as a
consequence of the previous estimates,  its contribution is at
least $\frac{q}{2}$.

\vskip2mm

Let $i_1 < i_2 < \dots < i_j$ be the indices in the third part. We
have

\begin{equation}  \label{theorem22} \sum_{g=1} ^j  m_{i_g} |S_{i_g}| \ge
\frac{q}{2} . \end{equation}

\noindent We are going to  consider two cases:

\vskip2mm

\noindent  (I) $2m_{i_j} > l$:   In this case  $|S_{i_j}| \ge
\frac{|A_2|}{40} \ge \frac{|A|}{80}$ and $\frac{l}{2} S_{i_j}$ is
a subset of $l^{\ast} A$. In view of the initial ideas presented
in the previous subsection, we set $A' = S_{i_j}$, $n' = 2n$ and
$l' =l/2$. Since $l^d |A| \ge 160 \times 2^d C n$

$$ l^{\prime d} |A'| \ge (l/2)^d \frac{|A|}{80} \ge \frac{1}{80 \times 2^d} l^d
|A| \ge 2Cn = Cn', $$

\noindent and

$$ l^{\prime d'} |A| = \Omega (l^{d'} |A|), $$

\noindent for any $1 \le d' \le d$. The last two estimates
guarantee that the triple $(A', n', l')$ is perfect and we are
done.

\vskip2mm

\noindent (II) $2m_{i_j} \le l$: In this case, we prove that
$l^{\ast}A$ contains an arithmetic progression of length $cl|A|$
(in other words, one can set the parameter $d'$ in Theorem
\ref{311} equal to  one). For any integer $a$ which is the sum of
$l-2m_{i_j}$ different elements in $A_1$ (the set we put aside at
the beginning of the proof), $a+ m_{i_j} S_{i_j}$ is a subset of
$l^{\ast} A$. On the other hand, as $|A|^2 \ge 80 n \log_2 n$,

$$m_{i_j} |S_{i_j}| \ge \frac{q}{4t} \ge \frac{|A|^2}{20 \log_2n }
\ge Cn $$

\noindent Theorem  \ref{111} implies that $m_{i_j} S_{i_{g}}$
contains an arithmetic progression of length

$$cm_{i_j} |S_{i_j}| \ge c \frac{q}{4t} \ge c \frac{|A|^2}{20 \log _2 n}  \ge  cl|A| $$

\noindent if $l \le |A| /20 \log_2 n$. The case when $l$ is larger
than $|A|/ 20\log_2 n$ requires an extra argument. Notice that by
the definition of the third partial sum and the assumption on
$|A|$

$$\frac{1}{2}m_{i_g} |S_{i_g}| \ge \frac{1}{2} \frac{q}{4t} \ge
\frac{|A|^2} {40 \log_2 n} \ge Cn. $$

\noindent Given this, we can  apply Theorem \ref{111} to
$\frac{1}{2}m_{i_g} S_{i_g}$ to obtain an arithmetic progression
of length $cm_{i_g} |S_{i_g}|$, for every index $g$ in the third
partial sum. To conclude, we use the following simple fact  to
glue these arithmetic progressions together

\begin {fact} \label{link} Any element in $\sum_{g=1}^j \frac{1}{2} m_{i_g} S_{i_g} $
can be represented by the sum of $m= \sum_{g=1}^j m_{i_g}$
different elements from $A_2$. \end{fact}

\vskip2mm \noindent {\bf Proof of Fact \ref{link}.} Greedy
algorithm. \hs

\vskip2mm It follows that  $\sum_{g=1}^j \frac{1}{2} m_{i_g}
S_{i_g} $ is a subset of $m^{\ast} A_2 $, with $m$ defined  as in
Fact \ref{link}. Finally, by applying Corollary \ref{sumtwo1}
iteratively one can show that $\sum_{g=1}^j \frac{1}{2} m_{i_g}
S_{i_g} $ contains an arithmetic progression of length

$$c \sum_{g=1}^j m_{i_g} |S_{i_g}| \ge c \frac{q}{2} \ge
c\frac{|A|^2}{10} \ge c \frac{l|A|}{5}. $$

\noindent Now we can add additional elements from $A_1$ to
$m^{\ast} A_2$ to obtain a subset of $l^{\ast} A$. \hs

\vskip2mm

\noindent This simple proof, unfortunately, cannot be repeated for
the case $|A| = o( \sqrt {n })$. However, the arguments presented
here will be useful later on.







\subsection {A structural lemma}

In view of the result in the previous subsection, we only have the
deal with the case $|A| = O( \sqrt {n \log n })$. Actually, this
upper bound on $|A|$ matters little, but it imposes a bound on $l$
that is critical. Notice that if $|A| = O( \sqrt {n \log n })$,
then in order to guarantee the assumption $l^d |A| \ge Cn$ of
Theorem \ref{311}, we must have

$$l =\Omega ( n^{1/2d- o(1)}) \gg \log^{10}_2  n. $$

\vskip2mm

In this subsection, we focus on those pairs $(l, A)$, where
$l^d|A|$ is close to $n$ (but not necessarily larger than $n$) and
$l$ is relatively large.  A key step in our proof is the
 following structural lemma, which asserts that if $l^{\ast}
A$ does not yield a proper GAP as claimed  by Theorem \ref{311},
then $A$ must contain a big subset which has a very rigid
structure.

\begin{lemma} \label{structural} For any positive constants $\nu$ and $d$ there
are  positive constants $\delta, \alpha$ and $d_1$  such that the
following holds. Let $A$ be  a subset of $[n]$, $l$ be a positive
integer and $ n \ge f(n) \ge 1$ be a function of $n$ such that
$$\max \{\log^{10} n, (40f(n)\log_2 n) ^{1/3d} \} \le l \le |A|/2$$
and $l^d |A| f(n) \ge n$. Then one of the following two statements
must hold

\begin{itemize}
\item $l^{\ast}A$ contains a proper GAP of rank $d'$ and volume
$\Omega ( l^{d'} |A|)$ for some $1 \le d' \le d$.

\item There is a subset $\tilde A$ of $A$ with cardinality at
least $\delta |A|$ which is contained in a GAP $P$ of rank $d_1$
and volume $O(|A| f(n)^{1+\nu} \log^{\alpha} n)$.

\end{itemize}

\end{lemma}

\noindent The function $f(n)$ can be seen as a {\it rigidity}
parameter. The closer $l^d |A|$ is to $n$, the more rigid is the
structure of $\tilde A$. With some extra work, the lower bound of
$l$ in the lemma can be improved: 10 can be replaced by any
constant larger than 1 and $1/3d$ can be replaced by any positive
constant. If we refine the result this way, the constants $\alpha,
\nu$ and $d_1$ will also depend on the new constants.

\vskip2mm  For the proof of Theorem \ref{311}, we only need the
special case when $f(n)=1$. We, however, choose to present Lemma
\ref{structural} in the above general form since it might be of
independent interest and the proof is not significantly harder
than that of the special case.

\vskip2mm With $f(n)=1$, Lemma \ref{structural} yields the
following corollary.

\begin{corollary} \label{structural1}  For any positive constant $d$ there
are positive constants $\delta, \alpha$ and $d_1$ such that the
following holds. Let $A$ be  a subset of $[n]$, $l$ be a positive
integer such that $l^d |A|  \ge Cn$. Then one of the following two
statements must hold

\begin{itemize}
\item $l^{\ast}A$ contains a proper GAP of rank $d'$ and volume
$\Omega ( l^{d'} |A|)$ for some $1 \le d' \le d$.

\item There is a subset $\tilde A$ of $A$ with cardinality at
least $\delta |A| $ which is contained in a GAP $P$ of  rank $d_1$
and volume $O(|A| \log^{\alpha} n)$.

\end{itemize}

\end{corollary}

\noindent Notice that the  set $\tilde A$ in Corollary
\ref{structural1} satisfies

$$l^d |{\tilde A}| \ge  l^d \delta |A| \ge \delta C n. $$

\noindent Since $\delta$ depends only on $d$, by increasing the
constant $C$ in Theorem \ref{311}, we can always assume that
$\delta C$ is sufficiently large. Thus, given Corollary
\ref{structural1}, it  suffices  to prove Theorem \ref{311} under
the additional condition that $A$ is a subset of density at least
$\frac{1}{\log^{\alpha} n}$ of a GAP of constant rank, where both
the rank and $\alpha$ are constants depending on $d$. We present
this proof in the next section. A reader who is eager to see this
proof  can delay the reading  of the rest of this section and jump
right to Section 7.

\vskip2mm

The rest of this section is devoted to the proof of Lemma
\ref{structural}. As this  proof is fairly long, we brake it into
four parts, each of which contains arguments of fairly different
nature. The main technical ingredient  of this proof is again a
tree argument, similar to what we used in the proof of Theorem
\ref{211}. However,  the algorithm here is more complicated than
the algorithm in Section 4, and the analysis is also more
challenging.

\vskip2mm

In order to set up the algorithm we first need to produce a large
amount of subsets of $A$ with a certain property. This will be
done in the next subsection. In subsection 7.10, we describe our
algorithm together with several simple observations. Subsection
7.12 is devoted to an inverse argument, which we use to derive the
desired  properties of $A$. This derivation is quite different
from and much more tricky than the one in Section 5. We wrap up
with the final subsection, subsection 7.14, which contains the
verification of an estimate claimed in subsection 7.12.

\subsection {Small sets with big sums}

The goal of this subsection is to show that any finite set $A$
contains a subset $B$ of small size ($O(\ln |A|)$) such that
$|l^{\ast} B|$ is large, where  $l = |B|/2$.

\begin{lemma} \label{311-0} Let $A$ be a finite set of real
numbers where $|A|$ is sufficiently large.  Then $A$ contains a
subset $B$ of at most $20 \log_2 |A|$ elements such that
$(\frac{|B|}{2})^{\ast} B$ has cardinality at least $|A|$.
\end{lemma}

{\bf \noindent Proof of Lemma \ref{311-0}.} We can assume, without
loss of generality, that $|A|$ is sufficiently large so that $|A|
\ge 100 \log_2 |A|$.  We choose the first two elements of $A$, say
$a_1, a_2$ arbitrarily. Once $a_1, \dots, a_{2i}$ have been
chosen, we next choose $a_{2i+1}$ and $ a_{2i+2} $ from $ A
\backslash \{a_1, \dots, a_{2i}\}$ such that

\begin{equation}   \label{1.1} |(i+1)^{\ast} \{a_1, \dots, a_{2i+1}, a_{2i+2}\}
|\ge 1.1 |i^{\ast} \{a_1, \dots, a_{2i}\}| \end{equation}

\noindent  (if there are many possible pairs, we choose an
arbitrary one). We stop at time $T$ when $|T^{\ast} \{a_1, \dots,
a_{2T} \}| \ge |A|$ and let $B =\{a_1, \dots, a_{2T} \}$. It is
clear that $|B| \le 2 \log_{1.1} |A| \le 20 \log_2 |A|$. The only
point we need to make now is to show that as far as $|i^{\ast}
\{a_1, \dots, a_{2i} \}| < |A|$, we can always find a pair
$(a_{2i+1}, a_{2i+2})$ to satisfy (\ref{1.1}). Assume  (for
contradiction) that we get stuck at the $i^{th}$ step and denote
by $S$ the sum set $i^{\ast} \{a_1, \dots, a_{2i} \}$. For any two
numbers $a , a' \in A \backslash \{a_1, \dots, a_{2i}\}$, $(a+ S )
\cup (a'+S)$ is a subset of $(i+1)^{\ast} \{a_1, \dots, a_{2i}, a,
a'\}$. So by the assumption we have

$$|(a+ S ) \cup (a'+S)| \le 1.1 |S|. $$

\noindent Since both $a+S$ and $a'+S$ have $|S|$ elements, it
follows that their intersection has at least $.9 |S|$ elements.
This implies that the equation $a'-a= x-y$ has at least $.9|S|$
solutions $(x,y)$ where $x \in S$ and $ y\in S$. Now let us fix
$a$ as the smallest element of $A\backslash \{a_1, \dots,
a_{2i}\}$ and choose $a'$ arbitrarily. There are $|A|-2i-1 \ge
.9|A|$ choices for $a'$, each of which  generates at least $.9|S|$
pairs $(x,y)$ where both $x$ and $y$ are elements of $S$. As all
$(x,y)$ pairs are different, we have that

$$.9|A| \times .9 |S| \le { |S| \choose 2}, $$

\noindent which implies that $|S| > |A|$, a contradiction. This
concludes the proof. \hs

\vskip2mm \noindent { \it Many good small sets.} Consider a set
$A$ as in Theorem \ref{311}. Applying  Lemma \ref{311-0}  to $A$
to obtain a small set $A_1$. Next, apply the lemma to $A\backslash
A_1$ to obtain a small  set $A_2$ and so on. Each time we add to
$A_i$ few ``dummy" elements to make its cardinality exactly $20
\log_2 |A|$. Stop when $A \backslash (\cup_{i=1}^m A_i)$ has less
than $2|A|/3$ elements for the first time. Without loss of
generality, we can assume that $20 \log_2 |A|$ is even and set
$l_0 =10 \log_2 |A|$. We have a collection $A_1, \dots, A_m$ of
disjoint subsets of $A$ with the following properties

\begin{itemize}

\item $|A_1| =\dots= |A_m| =   20 \log_2 |A| = 2l_0$.

\item $|l_0^{\ast} A_i| = \Big|\frac{|A_i|}{2}^{\ast} A_i \Big|
\ge (2/3-o(1))|A| > {|A|}/2 $.

\item $|A \backslash (\cup_{i=1}^m A_i)| = (2/3+o(1)) |A|$.

\end{itemize}

\noindent Here we assume that $\log_2 |A| = o(|A|)$ which explains
the error terms $o(1)$ in the last two properties. In the next
subsection, we consider an algorithm which uses the sets $A_i$ as
input.

\subsection {The algorithm}

\noindent Set  $B_i = l_0 ^{\ast} A_i$ for all $1 \le i \le m$. We
now give a description of our algorithm. This algorithm constructs
a subset of $l^{\ast} A$ in a particular way. We shall exploit the
fact that the cardinality of this subset is at most $|l^{\ast}A|
\le ln$ (since $l^{\ast}A$ itself is a subset of $[ln]$) in order
to derive information about $A$.

\vskip2mm

{\it \noindent The algorithm.}  To start, set $m_0=m$. Truncate
the set $B_i$'s so each of them has exactly $b_0= |A|/2$ elements.
Denote by $B^0_i$ the truncation of $B_i$. We start with the
sequence of sets $B^0_1, \dots, B^0_{m_0}$, each of which has
exactly $b_0$ elements. Without loss of generality, we may assume
that $m_0$ is a power of 4.  At the beginning, we call the
elements in $A^{[1]} =A \backslash (\cup_{i=1}^m A_i) $  {\it
available}.

\vskip2mm

A general step of the algorithm functions as follows. The input is
a sequence $B^t_1, \dots, B^t_{m_t}$ of sets of the same
cardinality $b_t$. Consider the sets $\cup_{h=1}^{K} (B^t_i +
B^t_j +x_h)$ where $1 \le i < j \le m_t$  and  $x_1, \dots, x_K$
are different available elements ($K$ is a large constant to be
specified later). Choose $i,j, x_1, \dots, x_K$ such that the
cardinality of $B'_1 = \cup_{h=1}^{K} (B^t_i + B^t_j +x_h)$ is
maximum (if there are many possibilities, choose an arbitrary
one). Remove $i$ and $j$ from the index set and the $x_i$'s from
the available set and repeat the operation to obtain $B'_2$ and so
on. We end up with a set sequence $B'_1, \dots, B'_{m_t/2 }$ where
$|B'_1| \ge \dots \ge |B'_{m_t/2 }|$.

Let  $m_{t+1} =m_t/4$ and  set $b_{t+1} = |B'_{m_{t+1}}|$.
Truncate $B'_i$'s ($i < m_{t+1}$) so that the remaining sets have
exactly $b_{t+1}$ elements each. Denote by $B^{t+1}_i$ the
remaining subset of $B'_i$. The sequence $B^{t+1}_1, \dots,
B^{t+1}_{m_{t+1}}$ is the output of the step.

\vskip2mm  If $m_{t+1} \ge 4$, then we continue with the next
step. Otherwise, the algorithm terminates.

\vskip2mm

\noindent Let us pause for a moment and make a series of
observations. All of these observations are easy to verify so we
omit their proofs.

\begin{itemize}
\item Define $l_{t+1} = 2l_t +1 $ for $t=0,1,2 \dots$. Then
$B^{t}_i$ is a subset of $l_t^{\ast} A$  for any plausible $t$ and
$i$.

\item As $A$ is a subset of $[n]$,  $B^t_i$ is a subset of $[l_t
n]$.

\item  For any plausible $t$, $ b_{t+1} \ge 2 b_t. $

\item   After each step, the length of the sequence shrinks by a
factor $4$.

\item  At the beginning we have $(2/3-o(1))$ available elements.
The number of elements $x_i$'s used in the whole algorithm is
$o(|A|)$, so at any step, there are always $(2/3-o(1))|A| $
available elements.
\end{itemize}

\noindent  Since $l \gg l_0 = O(\log_2 n)$, we can assume, without
loss of generality, that $l/l_0 $ is a power of two, $l/l_0 =
2^{s_2}$. Recall that $m_0=m \approx \frac{1}{3} \frac{|A|}{2l_0}$
($m$ is slightly larger than $\frac{1}{3} \frac{|A|}{2l_0}$) and
$|A| \ge 2l$. It follows that $l/l_0 \le 4 m_0$. As we assume
$m_0$ is a power of 4, $m_0 =4^{s_1}$, it follows that $2(s_1+1)
\ge s_2$.

\vskip2mm We set  $K = 2^{c_1d}$, where $c_1$ is a constant at
least $9$. We first claim that

\begin{equation}  \label{defK} (K/2) ^{s_2/2} > 40 l^d f(n) \log n. \end{equation}

\noindent Indeed, observe that

\begin{equation}  \label{K1} (K/2)^{s_1/2} \ge (2^{9d+1}/2)^{s_2/2} = 2^{9d
s_2/2}. \end{equation}

\noindent Recalling the definition of $s_2$, $2^{s_2}= l/l_0$. We
assume that $l \ge \log_2^{10}  n \ge l_0^9$, so $2^{s_2} \ge
l^{8/9} $. It follows that

$$ 2^{9d s_2/2} \ge l^{ 4 d} = l^d \times l^{3d}  \ge 40 l^d f(n)  \log n, $$

\noindent by the assumption on $l$.

\vskip2mm

We next prove the following fact.

\begin{fact} \label{fact311-1} There is an index $k \le s_2/2$ such that
$b_k \le K ^k b_0$.
\end{fact}

\vskip2mm

{\bf \noindent Proof of Fact \ref{fact311-1}.} By the second
observation we have that $b_k \le l_k n$ for any $k$. From the
definition of $l_t$ it is easy to prove (using induction) that

$$l_k \le 2^k l_0 + 2^k  \le 2^{k+1} l_0. $$

\noindent It follows that $b_k \le 2^{k+1} l_0n$ for any $k$.
Recall that $b_0 =|A|/2$ and $l_0 =10 \log_2 |A|$. If $b_k > K^k
b_0$, then we should have

$$ K^k |A|/2 = K^k b_0 < b_k  \le 2^{k+1} l_0 n \le 2^{k+1}
n \times (10 \log_2 |A|), $$

\noindent which implies

$$ (K/2)^k |A| < 40 n \log _2|A| \le 40 n \log_2 n. $$

\noindent On the other hand, (\ref{defK})  and the assumption $l^d
|A| f(n) \ge n$ of Lemma \ref{structural} together  imply

$$ (K/2)^{s_2/2} |A| \ge  40 l^d |A| f (n) \log_2 n \ge 40  n \log_2 n, $$

\noindent which is a contradiction. The proof is thus complete.
\hs









\subsection {The inverse argument}  Let $k$ be the first index where
$b_k \le K^k b_0$. This means $|B^{k}_{m_k} |\le K^k b_0$. By the
description of the algorithm

\begin{equation}  \label{lemma100} B^k_{m_k} =   \cup_{h= 1}^K (B^{k-1}_i +
B^{k-1}_j+x_h) \end{equation}

\noindent for some $i,j$ and $x_h$'s. Given (\ref{lemma100}), we
are going to exploit the bound  $|B^{k}_{m_k} |\le K^k b_0$ in
many ways. First, this bound and the definition of $k$ means that
$|B^{k-1}_i + B^{k-1}_j|$ is relatively small and so we can use
Freiman's theorem to derive some facts about the sets $B^{k-1}_i $
and $B^{k-1}_j$. Next, (\ref{lemma100}) and the bound on
$|B^{k}_{m_k}|$ imply that there should be a significant overlap
among the sets $(B^{k-1}_i + B^{k-1}_j+x_h)$'s. Thus, there should
be a correlation between the (available) elements $x_h$'s. This
correlation eventually leads us to a structural property of the
set of available elements. The set $\tilde A$ claimed in the lemma
will be a subset of this set.

\vskip2mm To start, notice that (\ref{lemma100}) implies

\begin{equation}  \label{lemma10} |B^k_{m_k}| \ge  |B^{k-1}_i + B^{k-1}_j| \end{equation}

\noindent where  $1 \le i <j \le m_{k-1}$ and both $B^{k-1}_i$ and
$B^{k-1}_j$ has cardinality $b_{k-1} \ge K^{k-1} b_0$. The
definition of $k$ then implies that $|B^k_{m_k}| \le K b_{k-1}$,
so

\begin{equation}  \label{lemma11}  |B^{k-1}_i + B^{k-1}_j| \le K |B^{k-1}_i|.
\end{equation}

\noindent Applying  Freiman's theorem to (\ref{lemma11}), we could
deduce that there is a generalized AP $R$ with constant rank
containing $B^{k-1}_i$ and  $\Vol (R)= O(|B^{k-1}_i|) =
O(b_{k-1})$.

We say that two elements $u$ and $v$ of $B^{k-1}_j$ are equivalent
if their difference is in $R-R$. If $u$ and $v$ are not equivalent
then the sets $u + B^{k-1}_i$ and $v+ B^{k-1}_i$ are disjoint,
since $B^{k-1}_i$ is a subset of $R$.  By (\ref{lemma11}), the
number of equivalent classes is at most $K$. Let us denote these
classes by $C_1, \dots, C_K$, where some of the $C_s$'s might be
empty. We have $B^{k-1} _i \subset R$ and $B^{k-1}_j \subset
\cup_{s=1}^K C_s$.

\vskip2mm Let us now take a close look at (\ref{lemma100}). The
assumption $|B^k_{m_k}| \le K|B^{k-1}_i|$ and  (\ref{lemma100})
imply that there must be a pair $s_1, s_2$ such that the
intersection $$(B^{k-1}_i + B^{k-1}_j + x_{s_1}) \cap (B^{k-1}_i +
B^{k-1}_j  + x_{s_2})$$ is not empty. Moreover, the set $ \{x_1,
\dots, x_K \}$ in (\ref{lemma100}) was chosen optimally. Thus, for
any set of $K$ available elements, there are two elements $x$ and
$y$ such that
$$(B^{k-1}_i + B^{k-1}_j + x) \cap (B^{k-1}_i + B^{k-1}_j  + y)$$ is
not empty. This implies

\begin{equation}  \label{lemma12} x - y \in  (B^{k-1}_i + B^{k-1}_j )-
(B^{k-1}_i + B^{k-1}_j ) \subset \cup_{1\le g,h \le K}
\Big((R+C_g)-(R+ C_h) \Big) . \end{equation}

Define a graph $G$ on the set of available elements as follows:
$x$ and $y$ are adjacent if and only if $x-y \in (B^{k-1}_i +
B^{k-1}_j )- (B^{k-1}_i + B^{k-1}_j )$. By the argument above, $G$
does not contain an independent set of size $K$, so  there should
be a vertex $x$ with degree at least $|V(G)|/K$. By
(\ref{lemma12}), there is a pair $(g,h)$ such that there are at
least $|V(G)|/K^3$ elements $y$ satisfying

\begin{equation}  \label{lemma13} x - y \in    (R+C_g)-(R+ C_h) . \end{equation}

Both $C_g$ and $C_h$ are subsets of translations of $R$; so the
set $\tilde A$ of the elements $y$ satisfying (\ref{lemma13}) is a
subset of a translation of $P= (R+R)-(R+R)$. Recall that at any
step, the number of available elements is $(1-o(1))|A_2|$, we have

\begin{equation}  \label{lemma14} |\tilde A| \ge (1-o(1))|A_2|/ K^3 =\Omega
(|A_2|). \end{equation}

Let us summarize what we have obtained here. We have found a
subset $\tilde A$ of $A$ of density at least $(2/3-o(1))/K^3
=\Omega (1)$ and a GAP $P$ which contains $\tilde A$.  In order to
complete the proof of the lemma, it remains to bound the volume of
$P$. We need to show that if the first statement of the lemma does
not hold, then

\begin{equation}  \label{volumeP} \Vol (P) = O(|A| f (n)^{1+\nu}  \log^{\alpha}
n). \end{equation}

\noindent At this point, we know that

\begin{equation}  \label{volumeP1} \Vol (P) =O(\Vol (R)) = O (b_k), \end{equation}

\noindent where $b_k \le K^k b_0 = K^k |A|$. Unfortunately, we
still do not know  much about $K^k$.  Our next task is to  prove
that if the first statement of the lemma does not hold, then

\begin{equation}  \label{volumeP2} K^k = O( f(n)^{1+ \nu}  \log^{\alpha} n),
\end{equation}

\noindent which implies (\ref{volumeP}).

\vskip2mm

In order to verify (\ref{volumeP2}), we need to exploit the
definition of the sets $B^k_{m_k}$ even more.  Notice that when we
define $B^k_{m_k}$, we choose $i$ and $j$ optimally.  On the other
hand, as $m_{k} =\frac{1}{4}  m_{k-1}$, for any remaining index
$i$, we have at least $m' =  m_{k-1}/2$ choices for $j$. This
means that there are $m'$ sets $B^{k-1}_{j_1}, \dots,
B^{k-1}_{j_{l_2}}$, all of the same cardinality $b_{k-1}$, such
that

\begin{equation}  \label{lemma15} |B^{k-1}_i + B^{k-1}_{j_{r}}| \le K b_{k-1}
\end{equation}

\noindent  for all $1 \le r \le m' $.

\vskip2mm  From now on, we work with the sets $B^{k-1}_{j_{r}}$,
$1\le r \le m'$. By considering  equivalent classes (as in the
paragraph following (\ref{lemma11})), we can show that for each
$r$, $B^{k-1}_{j_{r}}$ contains a subset $D_r$ which is a subset
of a translation of $R$ and $ |D_r| \ge |B^{k-1}_{j_{r}}| /K
=\Omega (\Vol (R))$.

\vskip2mm By Lemma  \ref{sum}, there is a constant $g$ such that
$D_1 + \dots + D_g$ contains a GAP $Q_1$ with cardinality at least
$\gamma \Vol (R)$ for some positive constant $\gamma$. Using the
next $g$ $D_i$'s, we can create $Q_2$ and so on. At the end, we
have $m^{''} =\lfloor m'/g \rfloor$ generalized AP $Q_1, \dots,
Q_{m^{''}}$. Each of these has rank $d_1 =rank (R)$ (this
parameter $d_1$ is irrelevant in the whole argument) and
cardinality at least $\gamma \Vol (R)$. Moreover, they are subsets
of translations of the GAP $R'= gR$ which also has volume $O(\Vol
(R))$.

\vskip2mm

Consider a GAP $Q_i$. Due to its large volume (compared to the
volume of $R'$), there are only $O(1)$ possibilities for its
difference set. Thus, there is a positive constant $\gamma_1$ such
that at least a $\gamma_1$ fraction of the $Q_i$'s has the same
difference set. Truncating if necessary, we can assume the
corresponding sides of these $Q_i$'s have the same length (the
truncation could decrease the volumes by at most a constant
factor). Since two GAP with the same difference sets and
corresponding sides having the same length are translations of
each other, we conclude that there is a GAP $Q$ (of rank $d_1$ and
cardinality at least $\gamma \Vol (R)$) and an integer $m^{'''} =
\Omega (m^{''})$ so that there are least $m^{'''}$ translations of
$Q$ among the $Q_i$'s. Without loss of generality, we can assume
that these translations are $Q_1, \dots, Q_{m^{'''}}$. Before
continuing, let us gather some facts about $Q_i$ and $m^{'''}$.

\begin{itemize}

\item $|Q|= |Q_i| = \Omega (\Vol (R)) = \Omega (b_k) = \Omega (K^k
b_0) \ge \beta K^k b_0,$ for some positive constant $\beta$.

\item $m^{'''} =\Omega (m^{''}) =\Omega (m^{'}) = \Omega (m_k) =
\Omega (m_0/4^k) \ge \mu m_0/4^k$, for some positive constant
$\mu$.

\end{itemize}

\vskip2mm

\noindent To proceed further, we need the following fact, whose
proof is delayed until the next subsection.

\begin{fact} \label{311-2} If $\Big( \frac{K}{2 \times 4^d}
\Big)^k \ge f(n) \log^{d+2} n$, then there is  $\bar l  \le l$
such that $(\bar l )^{\ast} A$ contains a proper GAP of rank $d'$
and volume $\Omega (l^{d'} |A|)$ for some $1\le d' \le d$.

\end{fact}

In order to have $l^{\ast} A$ instead of $(\bar l )^{\ast} A$ one
can do the usual ``reserving" trick. Prior to Fact \ref{link}),
put aside $l$ elements from $A$ for reserve. Repeat the whole
proof with the remaining set until Fact \ref{311-2}. Now, choose
$l-\bar l $ arbitrary elements from the reserved set and add their
sum to the set $(\bar l )^{\ast} A$ obtained in Fact \ref{311-2}.
The resulting set is a subset of $l^{\ast} A$ and it contains a
proper GAP as claimed in Theorem \ref{311}.

\vskip2mm Now we conclude the proof of the lemma via Fact
\ref{311-2}. If we assume that the first statement of the lemma
does not hold, then  by this fact we have that

$$\Big( \frac{K}{2 \times 4^d} \Big)^k < f(n)  \log^{d+2} n. $$

\noindent Recall that we set $K= 2^{c_1d}$ where $c_1$ is a
constant. By setting $c_1$ sufficiently large compared to $1/\nu$,
it follows that

$$K^k \le f(n)^{1+\nu} \log_2^{\alpha} n, $$

\noindent for some constant $\alpha= \alpha (\nu, d)$, proving
(\ref{volumeP2}). \hs

\subsection {Proof of Fact \ref{311-2}} To prove Fact \ref{311-2}, let us set $l^{'} = \ep \min
( \frac{l}{l_k}, m_k/2 )$, where $\ep$ is a sufficiently small
positive constant. Without loss of generality, we can assume that
$l'$ is an integer. The definition of $l'$ and the construction of
the $Q_i$'s imply that for a proper choice of $\ep$, $l'Q$ is a
translation of a subset of $(\bar l)^{\ast} A$ for some $\bar l\le
l$. Fact \ref{311-2} follows from Theorem \ref{111} and the
following

\begin{fact} \label{311-3} If $\Big( \frac{K}{2 \times 4^d} \Big)^k
\ge f(n) \log_2^{d+2} n$, then the following two inequalities hold

\begin{equation}  \label{311-31} (l')^d |Q| \gg l_k n  \end{equation}

\begin{equation}  \label{311-32} (l')^{d'} |Q| \ge l^{d'} |A|, 1\le d' \le d
\end{equation}

\noindent where $(l')^d |Q| \gg l_k n$ means that $\frac{(l')^d
|Q| } {l_k n}$ tends to infinity with $n$.

\end{fact}

\noindent We need to define $l'$ as above due to the following
reason. The tree might be too tall (having much more than $\log_2
(l/l_0)$ levels) or too short (having less than $\log_2(l/l_0)$
levels). In the first case we have to look at some immediate level
between the root and the leaves. This corresponds to the case $l'
= \ep (l/l_k)$. In the second case, we look at some level very
close to the root and this corresponds to the definition $l' =\ep
(m_k/2)$.

\vskip2mm

{\bf \noindent   Proof of Fact \ref{311-3}.} Consider an arbitrary
integer $d'$ between 1 and $d$. The definition of $l'$ naturally
leads to the following two cases:

\vskip2mm

{\bf \noindent Case 1.} $l/l_k \le m_k /2$. In this case $l'= \ep
(l/l_k)$. Recalling that there is a constant $\beta$ such that
$|Q| \ge \beta K^k b_0$ (see the paragraph preceding Fact
\ref{311-2}), we have that for any $d' \ge 1$

\begin{equation}  \label{311-33} (l')^{d'} |Q| \ge \ep^{d'} (\frac
{l}{l_k})^{d'} \times  \beta K^k b_0 =  \frac{\ep^{d'} \beta}{2}
l^{d'} |A| \frac{K^k} {l_k^{d'}}, \end{equation}

\noindent where in the last equation we use the fact that $b_0 =
|A|/2$. On the other hand, recall that $l_0 = 10 \log_2 |A|$, we
have

$$l_k \le 2^{k+1} l_0 = 20 \times 2^k \log_2 |A| \le 20 \times 2^k \log _2 n $$

\noindent  So, it follows from (\ref{311-33}) that for any $1 \le
d' \le d$

\begin{equation}  \label{311-34} (l')^{d'} |Q| \ge  \frac{ \ep^{d'}
\beta}{2\times 20^{d'} } l^{d'} |A| \frac{K^k} { 2^{kd'} \log_2
^{d'} n} \ge  \frac{\ep^d \beta}{2 \times 20^d} l^{d'} |A| \Big(
\frac{K}{2^d} \Big)^k \frac{1}{\log^d_2 n}, \end{equation}

\noindent where the second inequality follows from the assumption
that $d' \le d$. The assumption on $K$ in Fact \ref{311-3} implies
that

$$\Big( \frac{K}{2^d} \Big)^k \ge f(n) \log^{d+2}_2 n >
  (\frac{\ep^d \beta}{2 \times 20^d})^{-1} \log^d_2 n, $$

\noindent  so the right most formula in (\ref{311-34}) is larger
than $ l^{d'} |A|$, for any $1 \le d' \le d$. This proves the
second inequality in Fact \ref{311-3}. To verify the first
inequality, notice that (\ref{311-34}) implies

\begin{equation}  \label{311-35}  \frac{ (l')^{d} |Q|}{l_k} \ge  \frac{\ep^d
\beta}{2 \times 20^d} l^{d} |A| \Big( \frac{K}{2^d} \Big)^k
\frac{1}{l_k\log^d n} \ge \frac{\ep^d \beta}{2 \times 20^d}  \Big(
\frac{K}{2^d} \Big)^k \frac{l^d |A| }{l_k\log^d n}. \end{equation}

\noindent Since $l_k \le 20 \times 2^k \log_2 n$, it follows that

\begin{equation}  \label{311-36}  \frac{ (l')^{d} |Q|}{l_k}  \ge \frac{\ep^d
\beta}{40 \times 20^d}  \Big( \frac{K}{2^{d+1}} \Big)^k
\frac{l^d|A|}{\log^{d+1} n}. \end{equation}

The assumption on $K$ implies that $\Big( \frac{K}{2^{d+1}}
\Big)^k \ge f(n) \log^{d+2} n $, so the right most formula in
(\ref{311-34}) is at least

$$   \frac{\ep^d
\beta}{40 \times 20^d} \frac{l^d|A| f(n) \log^{d+2} n}{\log^{d+1}
n } \gg n , $$

\noindent due to the assumption $l^d |A| f(n) \ge n$ of Lemma
\ref{structural}. This verifies the first inequality and completes
the treatment of Case 1.

\vskip2mm

{\bf \noindent Case 2.} $l/l_k  > m_k /2$. In this case $l'= \ep
(m_k/2)$. Since $m_k = m_0/4^k$ and

$$m_0 \ge |A|/6l_0 = |A|/60 \log_2 n$$

\noindent  we have that

$$ l' \ge \frac{\ep m_0}{2 \times 4^k} = \frac{\ep |A|}{120 \times
4^k \log_2 n}. $$

\noindent So for any $1 \le d' \le d$

\barray  (l')^{d'} |Q| \ge \Big ( \frac{\ep|A|}{4^k\times 120 \log
n} \Big)^{d'}  \beta K^k \frac{|A|}{2} &=& \frac{\ep^{d'} \beta}{2
\times 120^{d'} } |A|^{d'+1} (\frac {K}{4^{d'}})^k \frac{1}{\log
^{d'} _2n} \\ &\ge& \frac{\ep^{d'} \beta}{2 \times 120^{d'} }
l^{d'} |A| (\frac {K}{4^{d}})^k \frac{1}{\log ^{d}_2 n}. \earray

\noindent Similar to the previous case, the assumption on $K$
guarantees that $(\frac {K}{4^{d}})^k \ge \log^{d+2}_2 n \gg
\log^d_2 n$ which implies that

$$ \frac{\ep^{d'} \beta}{4 \times 120^{d'} } l^{d'} |A| (\frac
{K}{4^{d}})^k \frac{1}{\log ^{d}_2 n} \gg l^{d'} |A|, $$

\noindent for any $1 \le d' \le d$, which proves the second
inequality in Fact \ref{311-3}. To verify the first inequality,
notice that

\begin{equation}  \label{lkk}  \frac{(l')^d |Q|}{l_k} \ge \frac{\ep^d \beta}{2
\times 120^{d} } l^{d} |A| (\frac {K}{4^{d}})^k \frac{1}{l_k
\log_2 ^{d} n}. \end{equation}

\noindent Similar to the pervious case, we use the estimate $l_k
\le 20 \times 2^k \log_2 n$. This and (\ref{lkk}) give

$$  \frac{(l')^d |Q|}{l_k} \ge \frac{\ep^d \beta}{40
\times 120^{d} } l^{d} |A| (\frac {K}{2 \times 4^{d}})^k \frac{1}{
\log ^{d+1} n}. $$

\noindent Here we need the full strength of the  assumption on
$K$: $(\frac {K}{2 \times 4^{d}})^k \ge \log^{d+2}_2 n$. From this
and the assumption that $l^d|A| f(n) \ge n$, it follows that

$$  \frac{(l')^d |Q|}{l_k} \ge \frac{\ep^d \beta}{40
\times 120^{d} } \frac{ l^{d} |A| f(n) \log^{d+2}_2 n} { \log_2
^{d+1} n} \ge \frac{\ep^d \beta}{40 \times 120^{d} } n \log _2 n
\gg n , $$

\noindent completing the proof. \hs

\section{Proof of Theorem \ref{311} (continued)}

Thanks to Corollary \ref{structural1}, from now on we can assume
that $A$ is a subset of a GAP $P$ of rank $d_1$ and volume at most
$|A| \log^{\alpha} n$, where where both $d_1$ and $\alpha$ are
constants depending on $d$. We  first use this structural property
to create a set $B$ whose elements have high multiplicity with
respect to $A$. The set $B$ is a candidate for the set $A'$ in a
perfect triplet that we desire. After having created $B$, the
remaining  (and also the hard) part of the proof is to show that
there is a sufficiently large $l' \le l/2$ such that each elements
of $l'B$ can be represented as a sum of $2l'$ distinct  elements
of $A$. This part requires a non-trivial extension of the tiling
argument used in our earlier paper \cite{szemvu2}. In order to
carry out this extension we need to prove some new properties of
proper GAPs.

\vskip2mm

This section is organized as follows. In subsection 8.1, we define
the set $B$ and derive several properties of this set. This
subsection also contains a proof of the theorem for the case when
$l$ is relatively small compared to $|A|$ (see Corollary
\ref{l1large}). The next subsection, subsection 8.3, is devoted to
the study of proper GAPs. The results of this subsection will be
used in  subsection 8.6 to prove further properties of the set
$B$. In subsection 8.7, we specify a  plan for constructing a
sumset $l'B$ as desired. This plan is executed in the next three
subsections, 8.8, 8.10 and 8.11. The final subsection, subsection
8.12, discusses  a common generalization of Theorem \ref{211} and
Theorem \ref{311}.

\subsection {Sets with high multiplicity} \label{B}

We are going to show that there is a large set every element of
which has high multiplicity with respect to $A$. Consider a
monotone sequence $m_1, m_2, \dots$ and let  $S_i$ be the set of
numbers with multiplicities between $m_i$ and $m_{i+1}$. A natural
way to find a large set with high multiplicity is  to set $m_i
=2^i$ and process as in subsection \ref{largeA}. Here, however, we
shall set the $m_i$'s somewhat differently, in order to serve a
purpose which will become clear later.

\vskip2mm

We define $m_i = \frac{|A|}{2^i i}$  for all $i= 1,2, \dots,
\log_2 |A|$ (observe that the sequence $m_i$ is decreasing). Let
$S_i$ be the set of those numbers whose multiplicities with
respect to $A$ is less than $m_i$ and at least $m_{i+1}$.  A
simple double counting shows

\begin{equation}   \label{multi1} \sum_{i=1}^{\log_2 |A|} m_i |S_i| \ge { {|A|
\choose 2}}. \end{equation}

Now we are going to make some use of the structure of $A$. Since
$A$ is a subset of a GAP $P$, $2A$ is a subset of $2P$. On the
other hand, as $P$ is a GAP of constant rank and volume $O(|A|
\log_2^{\alpha} n)$, so $2P$ is a GAP with the same rank and
volume $O(|A| \log_2^{\alpha} n)$. The set $S_i$ (for all $i$) is
a subset of $2A$, so it follows that

\begin{equation}  \label{multi2} |S_i| \le |2A| \le |2P| = O(\Vol (2P)) = O(|A|
\log^{\alpha}_2 n). \end{equation}

\noindent  By (\ref{multi2}), the sum of those $m_i |S_i|$ where
$m_i \le \frac{|A|} { \log^{\alpha +2}_2 n}$ is at most

\begin{equation}  \label{multi21} O \Big(  \frac{|A|} { \log^{\alpha+2}_2 n}
|A| \log^{\alpha}_2 n \Big) \times \log |A|  = O(( \frac{|A|} {
\log^{\alpha+2}_2 n} |A| \log^{\alpha +1}_2 n \Big) = o(|A|^2).
\end{equation}

\noindent This estimate allows us to omit these terms from the sum
in (\ref{multi1}) and so significantly reduce the number of terms
in the sum. Notice that for any $i > \log_2 \log_2^{\alpha +2} n
$, $m_i \le \frac{|A|} { \log_2^{\alpha +2} n}$, so we only have
to look at the small $i$'s, $i \le  \log_2 \log_2^{\alpha +2} n$.
From (\ref{multi1}) and (\ref{multi21}), we have

\begin{equation}  \label{multi3} \sum_{i=1}^{  \log_2 \log_2^{\alpha +2} n}
\frac{|A|}{i 2^i }  |S_i | = \sum_{i=1}^{  \log_2 \log_2^{\alpha
+2} n} m_i |S_i | \ge {|A| \choose 2} - o(|A|^2) =
(\frac{1}{2}-o(1)) |A|^2. \end{equation}

\noindent The fact that $\sum_{i=1}^{\infty} \frac{1}{i^2}
=\pi^2/6$ and (\ref{multi3}) imply that there should be an index
$1\le i \le \log \log^{\alpha +2} n$ so that

$$|S_i| \ge \frac{6}{\pi^2} \frac{2^i}{i} (\frac{1}{2} -o(1))  |A| > \frac{2^i}{4i} |A|. $$

\noindent Choose the smallest $i$ satisfying the above inequality
and rename the corresponding set $S_i$ to $B$. We are going to
work with $B$ in the rest of the proof.  We set  $l_1 = \frac{|A|}
{(i+1) 2^{i+1}}$. Since we shall use the letter $i$ as an index
later, let us set $t= 2^{i+1}$ to avoid confusion. Under this new
notation, $l_1 = \frac{|A|} { t \log_2 t }$, where $t = 2^{i+1}$
is at most $2^{ \log_2 \log_2^{\alpha+2} |A|} \le \log_2^{\alpha+2
} n$. By the definition of the $S_i$'s, every element of $B$ has
multiplicity at least $l_1$ with respect to $A$. This implies that
$kB$ is a subset of $(2k)^{\ast} A$ for any $k \le l_1$. Now let
us consider two cases:

\vskip2mm

{\bf \noindent Case 1:} $l \le 2l_1$. In this case, we  set
$A'=B$, $l'=l/2$ and $n'= 2n$ and follow the plan described in
subsection 7.2. It is easy to verify  that the triplet $(A',l',
n')$ is perfect. Thus we have the following corollary which proves
Theorem \ref{311} for the case $l$ is relatively small compared to
$|A|$.

\begin {corollary} \label{l1large}
For any fixed positive integer $d$ there are positive constants
$C, c$ and $\beta$ depending on $d$
 such that the following holds. For any positive integers $n $
and $l$ and any set  $A \subset [n]$ satisfying $l \le
\frac{|A|}{\log^{\beta} n}$ and $l^{d} |A| \ge C n$,  $l^{\ast} A$
contains a proper GAP of rank $d'$ and volume at least $cl^{d'}
|A|$, for some $1 \le d' \le d$.

\end{corollary}
In the remaining part of the paper, we consider the case $l \ge
2l_1$. Before going to the next subsection, let us summarize what
we have at this stage. We have created  a set $B \subset 2A
\subset [2n]$ where

\begin{itemize}

\item $B$ has at least $ \frac{|A|t } {4\log _2t}$ elements.

\item Each element of $B$ has multiplicity at least $l_1
 =\frac{|A|} {t \log _2t} $ with respect
to $A$.

\item $t \le \log_2 ^{\alpha +2} n$.
\end{itemize}

\subsection{ Proper GAPs revisited}

If $A$ and $2A$ is a subset of a normal GAP $Q$, it is tempting to
conclude that $A$ is a subset of $\frac{1}{2} Q$. A naive ``proof"
would go as follows: Assume that there is an element $x \in A
\backslash \frac{1}{2} Q$. Since $A \subset Q$, $x \in  Q
\backslash \frac{1}{2} Q $ and so $2x \in 2Q \backslash Q$. But
$2x \in 2A \subset Q$, a contradiction.

\vskip2mm

The trap is in the second sentence. Reasonable it sounds, the
statement  ``$x \in Q \backslash \frac{1}{2} Q $  implies $2x \in
2Q \backslash Q$" is not true.  It is not hard to work out an
example where $2x \in Q \cap 2Q$. We can, however, easily avoid
this subtlety. If we assume that $2Q$ is proper then $x \in Q
\backslash \frac{1}{2} Q $ indeed implies that  $2x \in 2Q
\backslash Q$. Thus we can conclude

\begin{fact} If  $A$ and $2A$ is a subset of a normal GAP $Q$ and $2Q$ is proper, then  $A$ is a subset of $\frac{1}{2} Q$.
\end{fact}

\vskip2mm

The above fact motivates the following lemma, which is the main
result of this subsection. We assume $Q$ is normal and its edges
are divisible by $l$, so $\frac{1}{l} Q$ can be defined.

\begin{lemma} \label{fraction} For any constants $d$ and $g$ there are
constants $\gamma$ and $ k$  such that the following holds. Let
$B$ be a finite set of integers, $l$  a positive integer and  $Q$
a (normal) proper GAP of rank $d$ satisfying

\begin{itemize}

\item The union of $g$ translations of $Q$ cover $lB$.

\item $kQ$ is proper.

\end{itemize}

\noindent Then there is a translation $B_1$ of $B$ such that $B_1
\cap \frac{1}{l}Q$  has at least $\gamma |B|$ elements.
\end{lemma}

{\noindent \bf Proof of Lemma \ref{fraction}.} We can assume,
without loss of generality, that $B$ contains 0. The normal GAP
$Q$ can be represented as  $Q= \{\sum_{i=1}^d x_i a_i| 0\le x_i
\le n_i\}$. If $lB$ is covered by $g$ translations of $Q$ then
$lB-lB$ is covered by $g_1= g^2$ translations of $P=Q-Q$, which
has the form $P= \{\sum_{i=1}^d x_i a_1| -n_i \le x_i \le n_i\}$.
 Let $P_1 = \frac{1}{2} P$ and $P_2 =\frac{1}{2} P_1$; it is
clear that $P_1$ is a translation of $Q$. Since $g_1$ translations
of $P$ cover $lB-lB$ and each translation of $P$ is the union of
$h^d$ translations of $P_1$, $lB-lB$ is covered by $2^d g_1$
translations of $P_1$. Furthermore, as each translation of $P_1$
is the union of $2^d$ translations of $P_2$, $lB-lB$ is covered by
$4^d g_1$ translations of $P_2$.

\vskip2mm Since $0 \in B$,  $lB-lB$ contains  $B$. By the pigeon
hole principle, there is a translation of $P_2$ containing at
least an $\frac{1}{4^d g_1}$ fraction of $B$. Equivalently, $P_2$
contains a set $B' \subset a + B$ where $|B'| \ge \gamma |B|$ and
$a$ is an integer. Setting $k= 2^{d+2} g_1$ and $h = 2^{d+1} g_1
+1$, we are going to show that $B'-B'$ is a subset of $\frac{h}{l}
P_1$. Since $B'-B'$ contains a subset of constant density of a
translation of $B$ and $P_1$ is a translation of $Q$, it follows
that  there is a translation  of $B$ which intersects
$\frac{h}{l}Q$ in  $\Omega(|B|)$ elements. This implies the claim
of the lemma since  $\frac{h}{l}Q$ is the union of $2^h =O(1)$
translations of $\frac{1}{l}Q$.

\vskip2mm In the rest of the proof, let us assume, for the sake of
a contradiction, that there is an element $x $ of $ B'-B'$ not
belonging to $\frac{h}{l} P_1$. Since  $B'-B'$ is a subset of
$P_2-P_2=P_1$, $x$ is an element of $P_1$.  Let $s_1$ be the
smallest positive integer such that $s_1 x \in 2P_1 \backslash
P_1$. Since both $2P_1$ and $P_1$ are proper, $s_1$ is at most
$l/h$.

\vskip2mm Recall that $B'$ is a subset of $a+B$. So, an element of
$B'$ has the form $a+b$ where $b \in B$. As $x \in B'-B'$,  $x
=b_1-b_2$ for some $b_1, b_2 \in B$. We set $y= s_1 x$ and
consider the sequence $y, 2y, 3y, \dots, \lfloor l/s_1 \rfloor$y.
As $s_1 \le l/h$, $\lfloor l/s_1 \rfloor \ge h > 2^{d+1} g_1$.
Each element of the above sequence has the form $r b_1-rb_2$ for
some $r \le l$. Since $0\in B$, these elements belong to $lB-lB$.
Let us now restrict ourself to the subsequence $$y, 2y, \dots,
(2^{d+1}g+1 +1)y. $$ Recall that $lB-lB$ is a subset of the union
of $ 2^d g_1$ translations of $P_1$.  The pigeon hole principle
implies that there should be a translation, say $a'+P_1$,
containing two elements $iy$ and $jy$ where $2 \le i-j \le 2^{d+1}
g_1$. The difference $(i-j)y$ is an element of $(a'+P_1)-(a'+P_1)=
P_1-P_1= 2P_1$.  Since $i-j < 2^{d+1} g_1=  k/2$, $2(i-j)P_1$ is
proper by the second assumption of the lemma. Moreover, $y$ is an
element of $2P_1 \backslash P_1$ so $(i-j)y$ is an element of
$2(i-j)P_1 \backslash (i-j)P_1$. This is a contradiction because
$(i-j)P_1$ contains $2P_1$ as $i-j \ge 2$. \hs

\subsection  {Properties of $B$} \label{propertyB}

Let us  consider the set $l_1B$. By the lower bounds on $l_1$ and
$|B|$ (see the last paragraph of subsection 7.1  we have

$$l_1^d |B| \ge \frac{|A|^{d+1}} { 4 t^{d-1} \log_2^{d+1} t } . $$

\noindent The assumptions $l^d|A| \ge Cn$ and $l \le |A|/2$ of
Theorem \ref{311} guarantee that $|A|^{d+1} \ge Cn$ and so

$$l_1^d |B| \ge \frac{C}{4 t^{d-1} \log_2^{d+1} t} n. $$

The factor $t^{d-1} \log_2^{d+1} t$ is the main source of our
troubles. If $t$ is a constant bounded by a function of $d$ (say
$e^{d^2}$), then by increasing the value of $C$  we can assume
that $ \frac{C}{4 t^{d-1} \log^{d+1} t}$ is sufficiently large and
so Theorem \ref{111} can by applied. However, $t$ can be as large
as a positive power of $\log _2 n$ and in general  cannot be
bounded by any function of $d$.

\vskip2mm In the remaining part of the proof, we assume that $t$
is very large compared to  $d$ (for all purposes, it is sufficient
to assume, say, $t \ge e^{e^{100d}}$). We are going to find a way
a play this assumption to our advantage (and through our arguments
one will see the reason for the somewhat artificial definition of
$m_i$'s). In the remaining part of this subsection, we use Lemma
\ref{fraction} to derive some properties of $B$ which are useful
for us.

\vskip2mm

Let us start with  the  usual ``doubling" trick. Set $B_0=B$ and
define $B_{i+1}= 2B_i$. We claim that at some stage we will be in
a position to apply Lemma \ref{fraction}.

\vskip2mm It is easy to show (using an argument similar to those
used in the proof of Theorem \ref{111}) that there is some  $s$
such that $2^s \ll l_1$  satisfying $|2B_s| \le (2^{d+2} -1)
|B_s|$. As usual, we let $s$ be the smallest number  with this
property. By Lemma \ref{bilu},  $B_s$ is a subset of a constant
number of translations of a GAP $P_0$ of rank $ d+1$ where $\Vol
(P_0) = O(|B_s|)$. Moreover, the proper filling lemma implies that
there is a constant $g_1$ so that  $g_1 B_s$ contains a proper GAP
$P_1$ of rank $d+1$ whose volume is $\Theta (|B_s|)$. The
differences of $P_1$ are constant  multiplies of the corresponding
differences of $P_0$, so $P_0$ is covered by a constant number of
translations of $P_1$. Therefore,  $B_s$ is covered by a constant
number of translations of the proper GAP $P_1$.

\vskip2mm

In order to apply Lemma \ref{fraction}, we also need the
assumption that there is a sufficiently large constant $k_1$ such
that  $k_1P_1$ is proper. Unfortunately, nothing guarantees the
existence of $k_1$. However, if we cannot find $k_1$, then we can
use our ``rank reduction" argument.  Set $k_1$ be a sufficiently
large constant and  consider the sequence $P_1, 2P_1, 4P_1,
\dots$. If for some $i \le \log_ 2k_1$, $2^iP_1$ fails to be
proper, then by the rank reduction argument, we can find a proper
GAP $P_2$ of rank strictly  less than the rank of $P_1$  such that
the following two properties hold

\begin{itemize}

\item There is a constant $g_2$ such that  $g_2 P_1$ contains
$P_2$.

\item A constant number of translations of $P_2$ cover $P_1$.

\end{itemize}

It follows that a constant number of translations of $P_2$ cover
$B_s$.  Now repeat the above argument with $P_2$. As the rank
decreases  each time, we should be done after a constant number of
steps. According to our arguments,  the final proper GAP (for
which the assumptions of Lemma \ref{fraction} are satisfied) still
has volume $\Omega (|B_s|)$. We call this final GAP $P'$.


\vskip2mm \noindent By applying Lemma \ref{fraction} to $P'$ we
obtain a few new properties of $B$

\begin{itemize}

\item For some $m = O(2^s) $, $mB$ contains a GAP $P'$ which has
volume at least

$$\Omega  (|B_s|) = \Omega  \Big((2^{d+2}-1)^s |B|\Big)  =\Omega  (2^{s(d+1)}
|B|). $$

\noindent Moreover, since $l_1 \gg 2^s$, $m \ll l_1$.

\item There is a subset $B'$ of $B$ such that $|B'| \ge \gamma
|B|$ and $B'$ is a subset of a GAP $P$ which is a translation of $
\frac{1}{m} P'$.

\end {itemize}

Since we are allowed to ignore constant factors, we  assume that
$B'=B$ for convenience.  Moreover, without loss of generality, we
could assume that $P'$ has symmetric form, namely, $P' =\{a_1x_1
+\dots a_{d_1} x_{d_1}| -n_i \le x_i \le n_i \}$.

\subsection  {A plan} \label{plan}

\noindent  Let us now give a rough discussion of our plan:

\begin{itemize}

\item  We are going to find  a set $\CT$ of $l_2$-tuples in $B$ (a
$k$-tuple is a set of $k$ not necessarily different elements) such
that the sum of the elements in any tuple is an element of
$(2l_2)^{\ast} A$, where $l_2 \gg l_1$ is a parameter to be
defined. Let $\CS$ be the collection of the sums of the tuples in
$\CT$. We create $\CT$ in a particular manner so that $\CS$ is
sufficiently dense in $l_2B$.

\item   We next prove that $ \CS + l_1B$ contains $l_2B$, relying
on  the fact that $\CS$ is dense in $l_2B$. This  way we obtain
the sum set $l_2B$ where $l_2$  is significantly larger than
$l_1$.

\item Since $\CS$ is a subset of $(2l_2)^{\ast}A$ and $l_1B$ is a
subset of $(2l_1)^{\ast} A$, $\CS + l_1B$ is a subset of
$(2l_2)^{\ast} A + (2l_1)^{\ast} A$. The obvious obstacle here is
that the same element of $A$ might be used twice, once in
$(2l_2)^{\ast} A$ and once in $(2l_1)^{\ast} A$. We overcome this
problem in  subsection 8.10 and show that $l_2B$ is in fact an
element of $(2l_2+ 2l_1)^{\ast} A$.

\end{itemize}

We call this plan a tiling operation as what it does is  to tile
many copies of $l_1B$ together to get a bigger set $l_2B$.

\vskip2mm

Would we be done after a successful implementation of this plan ?
Well, we would be  in a very good position if we can guarantee
that $l_2^d |B| \gg n$ (this inequality is necessary for an
application of Theorem \ref{111} to $l_2B$).   In the case  $d=1$,
we can do this  and the above plan was carried out successfully in
an earlier paper \cite{szemvu}. Unfortunately, there is  a serious
difference between  the two cases case $d=1$ and $d \ge 2$. For
$d=1$, the troublesome factor $t^{d-1} \log_2^{d+1} t$ is only
$\log_2^2 t$ and there is a way to set up $l_2$ so this
poly-logarithmic factor can be ignored. On the other hand, in the
general case $d \ge 2$, the troublesome factor is a polynomial in
$t$ (which is of a different order of magnitude) and even the
optimal value we could get for $l_2$ would not be enough  to kill
this factor.

\vskip2mm

We are going to resolve this problem by repeating the second step
of the plan  many times. Roughly speaking,  what we shall do is to
put many original tiles (copies of the set $l_1B$) together to get
a larger tile $l_2 B$. Next, we put many copies of $l_2 B$
together to get an even larger tile $l_3 B$ and so on. We repeat
the operation  until we get a sufficiently large tile $l_kB$ which
satisfies $l_{k}^d|B| \gg n$.

\vskip2mm There is a trade-off in this argument.  The repetitions
make the problem mentioned the last step of the above plan more
severe: Now the same element of $A$ might be used as many as $k$
times. Luckily, our treatment for this problem is not sensitive to
this modification as far as $k$ remains a constant, which is the
case.

\vskip2mm

Finally, let us go back to  address the first step: How can we
find $l_2$ elements of $B$ such that their sum can be represented
as the sum of $2l_2$ different elements of $A$ ? The main idea is
 as follows: An element of $B$ has multiplicity $l_1$ with
respect to $A$, so it gives us $l_1$ pairs of elements of $A$, all
have the same sum. Therefore, a set of $m$ different elements of
$B$ gives us $l_1 m$ different pairs. On the other hand, each
element in $A$ occurs in at most $|A|-1 < |A|$ pairs. Using the
greedy algorithm, we can find  at least $\frac{l_1m}{2|A|}$
mutually disjoint pairs. Thus,  for any $l_2 \le
\frac{l_1m}{2|A|}$, we have a collection of $l_2$ mutually
disjoint pairs. Clearly, the sum of the $l_2$ elements of $B$
corresponding to these pairs is an elements of $(2l_2)^{\ast} A$.

\vskip2mm

The critical feature of this step is how to choose the  set of $m$
elements of $B$. We discuss this issue  in   the next paragraph.

\subsection {The Tiling Operation: Start} \label{start}

Let us start with the execution of   the first step. Recall, from
the last paragraph of subsection \ref{propertyB}, that $B$ is a
subset of a proper GAP $P$ of constant rank $d_1$ (the value of
$d_1$ is irrelevant but we do know that $d_1\le d+1$). It is
easier for the reader to visualize the argument if he/she
identifies  $P$ with a $d_1$ dimensional box. Partition each edge
of $P$ into $T_1$ intervals of equal length, where $T_1$ is a
parameter to be determined. The products of these intervals
partition $P$ into $(T_1)^{d_1}$ identical small boxes. A small
box $Q$ is {\it dense} if the number of elements of $B$ in $Q$ is
at least $\frac{|B|}{2 (T_1)^{d_1}}$; $Q$ is {\it sparse}
otherwise. The sparse boxes contain at most half of the elements
of $B$, so at least half of the elements of $B$ should be
contained in dense boxes. Since constants like $1/2$ do not play
any significant role, we assume, for the sake of convenience, that
all elements of $B$ are contained in dense boxes.

\vskip2mm Let us recall that $|B| \ge \frac{|A|t}{4 \log _2t} $
and $l_1 = \frac{|A|}{t \log _2t}$. By throwing away dummy
elements, we can assume that $|B|$ is exactly $\frac{|A|t}{4
\log_2 t}$.

\vskip2mm

Consider a dense box $Q$, for each element $x \in B \cap Q$, $x$
has multiplicity $l_1$ with respect to $A$. We set the number $m$
in the last paragraph of the previous subsection to be $|B|/
2T_1^{d_1}$; as $Q$ is dense we are guaranteed to find this many
elements of $B$ in $Q$.  The argument in the above mentioned
paragraph shows that we can have at least

$$  \frac{l_1 |B|}{ 4 (T_1)^{d_1} |A|}, $$

\noindent disjoint pairs. For a technical reason, we do not set
$l_2 $ equal this value, but equal one-third of it:

$$l_2 =  \frac{l_1 |B|}{ 12 (T_1)^{d_1} |A|}. $$

\noindent For $x \in B$ let $N_x$ be the collection of pairs (in
$A$) summing up to $x$. We have proved

\begin{fact} \label{disjoint} For each dense box $ Q$, the union of
$N_x$'s ($x\in B \cap Q$) contains at least $3l_2$ mutually
disjoint pairs. \end{fact}

\noindent Substituting the values of  $l_1$ and $|B|$ into the
formula of $l_2$, we have

\begin{equation}  \label{l2}  l_2 = \frac{|A| } {48 (T_1)^{d_1} \log_2^2 t}.
\end{equation}

For each dense box $Q$, fix a collection  $N_Q$  of $3l_2$
disjoint pairs. For a pair $(a,b)$ in $N_Q$, the number $a+b$ is a
point of the box $Q$ ($a+b \in B \cap Q$). In the following, we
denote by $D_Q$ the collection of these points; $D_Q$ is a
multi-set as different pairs may have the same sum. Let $D$ be the
union of the $D_Q$'s.

\vskip2mm  Let us now take a closer look at the set $l_2B $. An
element $x$  of this set can be written as $x=x_1 +\dots +x_{l_2}
$, where $x_i$'s are not necessarily different elements of  $B$.
Moreover, we assumed that every element of $B$ is in some dense
box, so  each $x_i$ is in some dense box $Q$ (different $x_i$'s
may, of course, belong to different boxes). Fix a dense box $Q$;
for each $x_i \in Q$, we are going to replace it by some $y_i \in
D_Q$. Now comes a very important point. Since $|D_Q| \ge 3l_2$ for
any dense box $Q$, we can replace $x_1, \dots, x_{l_2}$ with
elements $y_1, \dots, y_{l_2}$ with the following property: There
are mutually disjoint pairs $(a_1, a_1'), \dots, (a_{l_2},
a_{l_2}')$, $a_i, a_i' \in A$, such that $a_i + a_i'= y_i$. To see
this, let us consider the following rule. For $x_1$, choose an
arbitrary pair $(a_1, a_1')$ from $D_{Q_1}$ where $Q_1$ is the
dense box containing $x_1$; set $y_1 = a_1 + a_1'$. Assume that
$(a_1,a_1'), \dots, (a_{i-1}, a_{i-1}')$ have been chosen.
Consider $x_i$ and the set $D_{Q_i}$ where $Q_i$ is the dense box
containing $x_i$. Delete from $D_{Q_i}$ every pair which has a
non-empty intersection with the chosen pairs. Since the pairs in
$D_{Q_i}$ are disjoint, any pair $(a_j,a_j')$ ($ 1\le j \le i-1$)
could intersect at most 2 pairs in $D_{Q_i}$ so we delete at most

$$ 2(i-1) \le 2(l_2 -1) < 2l_2 , $$

\noindent pairs from $D_{Q_i}$. But $D_{Q_i}$ contains $3l_2$
pairs so there are always some pairs left and we choose an
arbitrary one among these.

\vskip2mm

The disjointness of the chosen pairs guarantees that $y=y_1+\dots
+ y_{l_2}$ can be represented as a sum of exactly $2l_2$ different
elements from $A$. Let $\CT$ denote the collection of the tuples
$(y_1, \dots, y_{l_2})$ and $\CS$ be the collections of their
sums. Following the plan, we next  show that $\CS + l_1B$ contains
$l_2B$.

\vskip2mm  Consider  $x= x_1 +\dots + x_{l_2}$. Since $x_i \in B$
and $B \subset P$, each $x_i$'s is an element of the box $P$ and
can be viewed as a point in $\BBZ^{d_1}$, so we can view $x$ as a
vector in $\BBZ^{d_1}$. By replacing $x_i$ with $y_i$, we obtain
another vector  $y = \sum_{i=1}^{l_2} y_i$. We are going to find a
box $P_1$ centered at the origin so that $P_1$ is a subset of
$l_1B$ and  the difference $x-y= \sum_{i=1}^{l_2} (x_i -y_i)$ a
vector in $P_1$. The union of the copies of  such a $P_1$ centered
at the points of $\CS$ cover $l_2B$. As $P_1 \subset l_1B$, it
follows that $l_2B \subset \CS + l_1B$, as desired.

\vskip2mm

The key observation in what follows is that  $x_i-y_i$ is small
because they are in the same small box (this is the main reason
why we partition $P$ into many small boxes). Let us fix an edge of
$P$ and assume that its length is $s_1$. The absolute value of the
component of $x_i-y_i$ in the direction of this edge is at most
$s_1/T_1$. It follows that the corresponding component of $x-y$ is
at most $l_2s_1/ T_1$. We are going to choose $T_1$ and define
$P_1$ so that that this bound is at most half the length of the
corresponding edge of $P_1$ ($P_1$ is centered at the origin).
This would imply that $P_1$ contains the vector $x-y$.

\vskip2mm

Now we are going to define $P_1$.  The last paragraph of
subsection \ref{propertyB} tells us that  $mB$ contains a GAP
$P'=mP$, for some $m \ll l_1$.
 Thus $l_1B$ contains the box $ \frac{l_1}{m} P' =l_1P$. This is
 our box $P_1$. Observe that  $P_1$'s edge in the relevant direction has
length $s_1 l_1$. In order to guarantee that this length is at
least twice $l_2s_1/T_1$, we should set $T_1$ so that

\begin{equation}  \label{T1}  s_1 l_1 \ge  \frac{2 l_2s_1}{ T_1} =
\frac{s_1|A|}{24 T_1^{d_1+1} \log^2_2 t}. \end{equation}

\noindent To satisfy (\ref{T1}), it is sufficient to set

$$T_1 = \Big( \frac{ |A|}{24 l_1 \log_2^2 t} \Big) ^{1/(d_1+1)}  =
\Big( \frac{t}{24 \log_2 t }\Big)^{1/(d_1+1)}, $$

\noindent since $l_1 = |A|/ t \log_1 t$. For the sake of a cleaner
calculation, we set $T_1$ a little bit larger

$$T_1 = t ^{1/(d_1+1)}. $$

\noindent Substituting the above value of $T_1$ into the
definition of $l_2$ in (\ref{l2}), we obtain

\begin{equation}  \label{l21}   l_2 =  \frac{|A| } {48 (T_1)^{d_1} \log_2^2 t}
= \frac{|A|} {48 t^{d_1/(d_1+1)} \log_2^2 t} \ge \frac{|A|}
{t^{d_1/(d_1+1)} \log_2^2 t}. \end{equation}

\noindent This $l_2$ is still not large enough, namely,  $l_2^d
|B|$ could still be smaller than $n$. Indeed, the above lower
bound on $l_2$ only guarantees that

\begin{equation}  \label{122}  l_2^d |B| \ge \frac{|A|^{d}}{t^{dd_1/(d_1+1)}
\log^{2d} t} \times  \frac{|A|t}{8 \log_2 t } = \Theta
(\frac{|A|^{d+1}}{t^{dd_1/(d_1+1)-1} \log_2^{2d+1} t}) ,
\end{equation}

\noindent where the right hand side can be significantly smaller
than $n$ if $|A| = O (n^{1/(d+1)})$ and $dd_1/(d_1+1)-1
>0$. Our plan is to  increase the value of $l_2$ by repeated
tiling.

\vskip2mm

To conclude this subsection, let us discuss the problem that the
same element of $A$ might appear twice in a representation of an
element of $l_2B$. Observe that $l_2B$ is subset of $(2l_1 +
2l_2)^{\ast}A$ and thus any element of $l_2B$ is a sum of $2l_1 +
2l_2$ elements of $A$. However, as we already pointed out, an
element of $A$ can appear twice, once in $(2l_1)^{\ast}A$ and once
in $(2l_1)^{\ast}A$. This problem can be resolved by the so-called
{\it cloning} trick, introduced in \cite{szemvu2}.

\subsection {The cloning argument} \label{clone}
At the very beginning of the whole proof, we split the set $A$
into two sets $A'$ and $A^{''}$ in such a way that $|A'| \approx
|A^{''}|$ and any number $x$ which has high multiplicity with
respect to $A'$ should have almost the same multiplicity with
respect to $A^{''}$. Next, we continue with $A'$ and keep $A^{''}$
for reserve. Repeat the whole proof with $A'$ playing the role of
$A$ until the previous paragraph. We call the set of elements with
high multiplicity (with respect to $A'$) $B'$ instead of $B$. Now
doing the same with $A^{''}$ we obtain a set $B^{''}$.

\vskip2mm

The key point now is that with  a proper splitting, the two sets
$B'$ and $B^{''}$ are exactly the same. So when we look at $l_2
B'$ as a subset of $\CS + l_1 B'$, we can think of an element of
$\CS$ as a sum of $l_2$ elements from $B^{''}$, rather than from
$B'$. Therefore, when we replace each element from $B'$ and
$B^{''}$ by the sum of two elements from $A$, the elements used
for $\CS$ come from $A^{''}$ and the elements used for $l_1B'$
come from $A'$ and this guarantees that no element of $A$ is used
twice. \vskip2mm

A random splitting provides the sets $A'$ and $A^{''}$ as
required. For each element of $A$ throw a fair coin. If head, we
put it into $A'$, otherwise it goes to $A^{''}$. If a number $x$
has multiplicity $m_x \gg \log n$ with respect to $A$, then
standard large deviation inequalities (such as Chernoff's) tell us
that with probability at least $1-n^{-2}$, $x$ has multiplicities
$\frac{m_x}{4} \pm 10 \sqrt {m_x \log n} = (1+o(1))\frac{m_x}{4}$
with respect to both $A'$ and $A^{''}$. Since there are only at
most $2n$ numbers $x$ to consider, with probability close to 1,
every $x$ with multiplicity $\gg \log n$ has approximately the
same multiplicities in $A'$ and $A^{''}$.

\vskip2mm  When we create the set $S_i$  (which we later rename to
$B$) in subsection \ref{B}, any element $x$ in $S_i$ has
multiplicity $m(x)$ at least $\frac{|A|}{2^{i+1} (i+1)} \gg \log
n$ with respect to $A$. So $x$ will have multiplicity roughly
$m(x)/4$ with respect to both $A'$ and $A^{''}$. Thus one can
expect that $x$ will appear in both $B'$ and $B^{''}$. The only
case we may have to worry about is when $m(x)$ is very close to a
threshold (say $m_i$) and then (because  the error terms can go
either way) $x$ might be in $B'$ but not in $B^{''}$ (or vice
versa). This problem is easy to deal with, we just force this $x$
to be in both $B'$ and $B^{''}$ (of course, forcing $x$ might
decrease $l_1$ slightly (by a factor $.9$, say) but this does not
influence anything).

\subsection {The tiling operation: Finish} \label{finish}

\noindent We repeat the tiling operation in  subsection
\ref{start} with new parameters. Now $P$ is cut into $T_2^{d_1}$
boxes, where $T_2$ is a parameter to be chosen. Instead of
(\ref{l2}), we define

\begin{equation}  \label{l3} l_3 =  \frac{|A|} {48 (T_2)^{d_1} \log_2 ^2 t }
\end{equation}

\noindent Here is our key point: in order to obtain $l_3B$, we now
add $\CS$ with $l_2B$, instead of with $l_1B$ as in  subsection
\ref{start}. This means that instead of $P_1$ we can use the
larger box $P_2= \frac{l_2}{l_1} P_1$. As an analogue of
(\ref{T1}), the condition we need on $T_2$ is

\begin{equation}  \label{T2}  s_1 l_2 \ge \frac{2 l_3 s_1}{ T_2} \end{equation}

\noindent  Notice that in the left hand side of (\ref{T2}) we have
$l_2$ instead of $l_1$. The fact that $l_2 \gg l_1$ allows us to
set $T_2$ much smaller than $T_1$. Consequently, $l_3$ becomes
significantly larger than $l_2$. Repeating this results in a
sequence $l_1 < l_2 < l_3 < l_4 < \dots$, where for some constant
$k$, $l_k$ will be sufficiently large.

\vskip2mm

Now let us present some computation. The derivation of $T_2$ from
(\ref{T2}) is similar to that of $T_1$ from (\ref{T1}). It is
sufficient to set

$$ T_2 =  \Big( \frac{ |A|}{24 l_2 \log_2^2 t} \Big) ^{1/(d_1+1)}  $$

\noindent in order to satisfy (\ref{T2}).  Since $l_2 \ge
\frac{|A|}{t^{d_1/d_1+1} \log^2_2 t}$,

$$ \Big( \frac{ |A|}{24 l_2 \log_2^2 t} \Big) ^{1/(d_1+1)}  \le
\Big( \frac{ t^{d_1/(d_1+1)}}{24 }\Big)^{1/(d_1+1)}, $$

\noindent so we can set $T_2 =\Big( \frac{t^{d_1/(d_1+1)} t}{24
}\Big)^{1/(d_1+1)}$.  Again, for convenience,  we set $T_2$ a bit
larger

$$T_2 = t^{d_1/(d_1+1)^2 }, $$

\noindent which implies

\begin{equation}  \label{l3-1}  l_3 =  \frac{|A| } {48 (T_2)^{d_1} \log^2_2  t}
\ge \frac{|A|} {t^{d_1^2/(d_1+1)^2} \log^2_2 t} . \end{equation}

\noindent By induction, we can show

\begin{equation}  \label{lk}  l_k \ge \frac{|A| } {t^{d_1^k/(d_1+1)^k }
\log^2_2 t}. \end{equation}

\noindent By choosing $k$ sufficiently large (say, $k = 2(d_1+1)
\log (d+1)$),  we have (using the fact that $t$ is much larger
than $d$)

$$l_k \ge \frac{|A|} {t^{1/2(d+1)}  \log^2 _2 t} \ge  \frac{|A|}
{t^{1/2d}}. $$

\noindent $l_k$ is now sufficiently large, namely, it satisfies
the critical inequality  $l_k^d |B| \gg n$ (one can easily check
this by substituting $|B| =\frac{|A|t}{4 \log t}$). This
inequality provides the necessary condition we need to  apply
Theorem \ref{111} to the set $l_k B$.

\vskip2mm

Our proof shows that $l_k B$ is a subset of $(2l_k)^{\ast} A +
(2l_{k-1})^{\ast} A + \dots + (2l_1)^{\ast} A$. In this sum an
element of $A$ might be used $k$ times. This problem can be
handled using the cloning argument exactly as before, with the
only formal modification that instead of splitting $A$ into two
subsets, we split it into $k$ subsets.

\vskip2mm

To be completely done, there is one last issue we need to discuss
and that is  the magnitude of the sum $l_1+ \dots+ l_k$.

\vskip2mm

As we have shown (with the aid of  cloning), the set $l_k B$ is a
subset of

$$(2l_1 + \dots + 2l_k)^{\ast} A = {\tilde l}^{\ast} A,$$

\noindent  where $\tilde l= 2l_1 + \dots + 2l_k$.
 We need to compare $\tilde l$ with $l$ and naturally
there are two cases. If ${\tilde l} \le l$, then we set $A'=B$,
$l'= \tilde l$ and $n'=2n$. In this case, we have

\barray  (l')^{d'} |A'| &\ge& l_k^{d'} |B|  \\ &\ge&
(\frac{|A|}{t^{1/2d}})^{d'} |B| \\&\ge&
(\frac{|A|}{t^{1/2d}})^{d'}\times \frac{|A|t}{4 \log_2 t} \\ &\ge&
{|A|^{d'+1} t^{1- (d'+1)/2d}}  \,\, \hbox{(as}
\,\, t \,\, \hbox{much larger than} \,\, \log_2 t) \\
&\ge& |A|^{d'+1} \\&\ge& l^{d'} |A|, \earray

\noindent for every $1 \le d' \le d$. This guarantees that the
triple $(A',l',n')$ is perfect.  \vskip2mm

In the remaining case when $\tilde l > l$, there is an index $i <
k$ such that  $$2l_1 + \dots 2l_i \le l < 2l_1 + \dots + 2l_{i+1}.
$$ We now modify  the tiling operation a little bit. First of all,
it is clear that we do not have to proceed beyond the $i$th tiling
so we make this tiling our last. Moreover, in this last tiling  we
shall not use the whole set $l_iB$ as a tile, but only a fraction
of it, say $l_i'B$ for some $l_i' < l_i$ (as we mentioned many
times, our arguments are invariant with respect translations so we
can assume that $l_i'B$ is a subset of $l_iB$). As the result, we
obtain  a set $l_{i+1}'B$ instead of $l_{i+1}B$, for some
$l_{i+1}'\le l_{i+1}$. The set $l_{i+1}'B$ is a subset of $(2l_1
+\dots + 2l_i + 2l_{i+1}') ^{\ast} A$ where, with a proper choose
of $l_i'$, we can guarantee that

$$  l/2 \le (2l_1
+\dots + 2l_i + 2l_{i+1}')  \le l. $$

\noindent Now we can set $A'=B$, $l'= l_{i+1}'$, $n'=2n$ and
conclude the proof as discussed in subsection \ref{perfect}. \hs

\subsection {A common generalization of Theorems \ref{211} and
\ref{311}}

In this subsection, we present a common generalization of Theorems
\ref{211} and \ref{311}. Let us first remind the reader of the
sumsets studied in these two theorems.  In Theorem \ref{211}, we
consider a sum of different sets $A_1, \dots, A_l$, but allow  the
same number to appear many times in a representation (the same
number may occur in several $A_i$'s). On the other hand,  in
Theorem \ref{311} we have only one set $A$ in the sum, but with
the restriction that the summands of a representation must be
different. For a common generalization of these theorems, we
consider a sum which involves different elements of different
sets. Let $A_1, \dots, A_l$ be sets of integers, we define $A_1
\stackrel{\ast}{+} A_2 \stackrel{\ast}{+} \dots \stackrel{\ast}{+}
A_l$ as the collection of all numbers which can be represented as
a sum of $l$ different numbers $a_1 \in A_1, \dots, a_l \in A_l$.
Formally speaking

$$ A_1 \stackrel{\ast}{+} A_2 \stackrel{\ast}{+}  \dots
\stackrel{\ast}{+} A_l = \{a_1 + \dots +a_l| a_i \in A_i, a_i \neq
a_j \,\, \hbox{for} \,\, 1\le i < j \le l \}. $$

We refer to $ A_1 \stackrel{\ast}{+} A_2$ as the star sum of $A_1$
and $A_2$.

\begin{theorem}  \label{411} For any fixed positive integer $d$
there are positive constants $C$ and $c$ depending on $d$
 such that the following holds. Let $A_1, \dots, A_l$ be
 subsets of size $|A|$ of $[n]$ where
 $l$ and $|A|$ satisfy $l^d |A| \ge Cn$.
Then $ A_1 \stackrel{\ast}{+} A_2 \stackrel{\ast}{+}  \dots
\stackrel{\ast}{+} A_l$ contains a GAP of rank $d'$ and volume at
least $cl^{d'} |A|$, for some integer $1 \le d' \le d$.
\end{theorem}

About the proof, one's first impression would be that  one can
prove Theorem \ref{411} using Theorem \ref{311} the same way one
proved Theorem \ref{211} using Theorem \ref{111}. This, however,
is not possible due to a subtle problem involving  star sums.
While it is clear that the (set) equality

 $$(A_1+ A_2) + (A_3+A_4)= A_1 + A_2 + A_3 +
A_4$$

\noindent  is true, its star sum counterpart

$$(A_1 \stackrel{\ast}{+}  A_2) \stackrel{\ast}{+} (A_3
\stackrel{\ast}{+} A_4)= A_1 \stackrel{\ast}{+} A_2
\stackrel{\ast}{+} A_3 \stackrel{\ast}{+} A_4 $$

\noindent is false.

\vskip2mm So far, the only way (we know) to verify  Theorem
\ref{411} is to repeat the proof of Theorem \ref{311} with
appropriate modifications. This is a tedious task, but no
essential new arguments are required, and  we thus omit the
details. Let us, however, present the variant of a step in the
proof of Theorem \ref{311}, Lemma \ref{311-0}, in order to give
the reader an idea about the kind of modifications one needs to
carry out.

\begin{lemma} \label{411-0} Let $A_i$, $1\le i \le 20 \log_2 |A|,$
be  finite sets of real numbers with the same cardinality $|A|$,
where $|A|$ is sufficiently large. Then there is an integer  $1
\le T \le 10 \log_2 |A|$ and elements $a_1 \in A_1, a_2 \in A_2,
\dots, a_{2T} \in A_{2T}$ such that all $a_i$'s are different and
the set $B=\{a_1, \dots, a_{2T} \}$ satisfies

$$|T^{\ast} B| \ge |A|. $$
\end{lemma}

{\bf \noindent Proof of Lemma \ref{411-0}.} We assume that $|A|$
is sufficiently large so that $|A| \ge 100 \log_2 |A|$. We choose
$a_1 $ and $a_2 $ from $ A_1$ and $A_2$, respectively, with the
only condition that $a_1 \neq a_2$. Once $a_1, \dots, a_{2i}$ have
been chosen, we next choose $a_{2i+1}$ and $a_{2i+2}$ from
$A_{2i+1}\backslash \{a_1, \dots, a_{2i}\}$ and
$A_{2i+2}\backslash \{a_1, \dots, a_{2i}\}$ so that $a_{2i+1} \neq
a_{2i+2}$ and

\begin{equation}   \label{4.1} |(i+1)^{\ast} \{a_1, \dots, a_{2i+1}, a_{2i+2}\}
|\ge 1.1 |i^{\ast} \{a_1, \dots, a_{2i}\}| \end{equation}

\noindent  (if there are many possible pairs, we choose an
arbitrary one). We stop at time $T$ when $|T^{\ast} \{a_1, \dots,
a_{2T} \}| \ge |A|$ and set  $B =\{a_1, \dots, a_{2T} \}$. It is
clear that $|B| \le 2 \log_{1.1} |A| \le 20 \log_2 |A|$. The only
point we need to make  is to show that as long as $|i^{\ast}
\{a_1, \dots, a_{2i} \}| < |A|$, we can always find a pair
$(a_{2i+1}, a_{2i+2})$ to satisfy (\ref{4.1}). Assume (for a
contradiction) that we get stuck at the $i^{th}$ step and denote
by $S$ the sum set $i^{\ast} \{a_1, \dots, a_{2i} \}$. For any two
numbers $a \in A_{2i+1} \backslash \{a_1, \dots, a_{2i}\}$, $a'
\in A_{2i+2} \backslash \{a_1, \dots, a_{2i}\}$ the union $(a+ S )
\cup (a'+S)$ is a subset of $(i+1)^{\ast} \{a_1, \dots, a_{2i}, a,
a'\}$. So by the assumption we have

$$|(a+ S ) \cup (a'+S)| \le 1.1 |S|. $$

\noindent Since both $a+S$ and $a'+S$ have $|S|$ elements, it
follows that their intersection has at least $.9 |S|$ elements.
This implies that the equation $a'-a= x-y$ has at least $.9|S|$
solutions $(x,y)$ where $x \in S$ and $ y\in S$. Now let us fix
$a$ as the smallest element of $A_{2i+1}\backslash \{a_1, \dots,
a_{2i}\}$ and choose $a'$ arbitrarily from $A_{2i+2} \backslash
\{a_1, \dots, a_{2i}, a\}$ (we exclude $a'$ from $A_{2i+2}$ so we
are guaranteed that $a \neq a'$). There are at least $|A|-2i-1 \ge
.9|A|$ choices for $a'$, each of which generates at least $.9|S|$
pairs $(x,y)$ where both $x$ and $y$ are elements of $S$. As all
$(x,y)$ pairs are different, we have that

$$.9|A| \times .9 |S| \le { |S| \choose 2}, $$

\noindent which implies that $|S| > |A|$, a contradiction. This
concludes the proof. \hs

\section {Erd\"os' conjecture on complete sequences}

In 1962, Erd\"os  introduced the following notion which has later
become quite popular: An infinite set $A$ of positive integers is
{\it complete} if every sufficiently large positive integer can be
represented as a sum of different elements of $A$ (see Section 6
of \cite{EG} or Section 4.3 of \cite{PS} for surveys about
completeness). For instance, Vinogradov's result (mentioned in the
Overview) implies that the set of primes is complete. On the other
hand, there is a big difference between the study of complete
sequences and the study of classical problems of Vinogradov-Waring
type. For completeness, we do not require the number of summands
in a representation to be the same. This relaxation leads to a
quite different kind of results. For problems of Vinogradov-Waring
type (where the number of summands is fixed), one usually requires
a very precise description  of the sequence (the set of primes or
the set of squares, say). For problems concerning complete
sequences, it has turned out there is much more flexibility.

\vskip2mm What would be the first condition for a sequence to be
complete ? Well, density must be the answer, as one cannot hope to
represent every positive integer with a very sparse sequence. But
one would also notice instantly that density itself would not be
enough: The set of even numbers has very high density, but is
clearly not complete. This shows that one should also consider a
condition involving modularity.

\vskip2mm In number theory it happens quite frequently that the
obvious necessary conditions are also sufficient. In 1962, Erd\"os
made the following conjecture

\begin{conj} \label{erd1} There is a constant $c$ such that the
following holds.  Any increasing sequence $A=\{a_1 < a_2 < a_ 3 <
\dots \}$ satisfying

(a) $A(n) \ge cn^{1/2}$

(b)  $S_A$ contains an element  of every infinite arithmetic
progression,

\vskip2mm

\noindent is complete.
\end{conj}

Here and later $A(n)$ denotes the number of elements of $A$ not
exceeding $n$. The bound on $A(n)$ is best possible, up to the
constant factor $c$, as shown by Cassels \cite{Cass}.

\vskip2mm  Erd\"os \cite{erd} proved that the statement of the
conjecture holds if one replaces (a) by a stronger condition that
$A(n) \ge cn^{(\sqrt {5}-1)/2}$. An important step was later made
by Folkman \cite{folk}, who improved Erd\"os' result by showing
that $A(n) \ge cn^{ 1/2 +\ep}$ is sufficient, for any positive
constant $\ep$. The first and simpler part in Folkman's proof is
to remove the condition (b). He showed that any sequence
satisfying (b) could be partitioned into two subsequences with the
same density, one of which still satisfies (b). In the next and
critical  step, Folkman shows that if $A$ is a sequence with
density at least $n^{ 1/2 +\ep}$ then $S_A$ contains an infinite
arithmetic progression (in other words, $A$ is subcomplete). His
result follows immediately from these two steps. Folkman's proof,
naturally, led him to the following conjecture, which is perhaps
even more to the point than Conjecture \ref{erd1}

\begin{conj} \label{folk1}  There is a constant $c$ such that the
following holds.  Any increasing sequence $A=\{a_1 < a_2 < a_ 3 <
\dots \}$ satisfying $A(n) \ge cn^{1/2}$ is subcomplete.
\end{conj}

\vskip2mm  Folkman's result has further been strengthened recently
by Hegyv\'ari \cite{Heg} and \L uczak and Schoen \cite{LS}, who
(independently) reduced the density $n^{1/2+\ep}$ to $cn^{1/2}
\log^{1/2} n$, using the result of  S\'ark\"ozy (see Section 3).

\vskip2mm In a previous paper \cite{szemvu2}, we proved Conjecture
\ref{folk1}. However, we decide to discuss this problem here for
pedagogical reasons. It would be more useful for the reader  to
consider this problem together with Conjecture \ref{folk} and
under  the general sufficient condition proved in Section 6. As a
matter of fact, given this sufficient condition, it is now very
simple to prove Conjecture \ref{folk1}. The only modification one
needs to make is to replace Lemma \ref{folkman1} by the following

\begin{lemma} \label{folkman1-1} There is a constant $C$ such that the following holds.
 If $A$ is a set of different positive
integers between 1 and $n$ and $|A| \ge C \sqrt n$, then $S_A$
contains an arithmetic progression of length $n$. \end{lemma}

\noindent The rest of the proof is the same.

\begin{theorem}  \label{folk1-10}  There is a constant $c$ such that the
following holds.  Any increasing sequence $A=\{a_1 < a_2 < a_ 3 <
\dots \}$ satisfying $A(n) \ge cn^{1/2}$ is subcomplete.
\end{theorem}

\vskip2mm Let us conclude with  a comment on Conjecture
\ref{folk1} and Conjecture \ref{folk}. These conjectures look
quite similar, which comes as no surprise as they appeared in the
same paper. The interesting point here is that the proof of
Conjecture \ref{folk} requires only Theorem \ref{211}, which is an
easy application of Theorem \ref{111}, but the proof of Conjecture
\ref{folk1} requires the much harder Theorem \ref{311}. On the
other hand, prior to  our study, Conjecture \ref{folk} seemed
harder to attack and less partial results were known.

\vskip2mm

{\bf \noindent Remark.} We have recently been informed by Lev
(private communication) that Chen \cite{Chen} also proved Theorem
\ref{folk1-10}, using a different method.

\section {Arithmetic progressions in finite fields}

In this section we assume that $n$ is a prime. We are going to
extends our previous theorems to arithmetic progressions modulo
$n$. The quantitative statements in these theorems will change
slightly, but the proofs remain essentially  the same. We first
establish the results and then describe an application.

\subsection {Results} In order to show why we need a modification in the
statements of the theorems,
 let us
consider the proof of Theorem \ref{111}. At one point in the proof
(see the paragraph following (\ref{ln})), we used the fact that
$lA$ is a subset of the interval $[ln]$ and thus has cardinality
at most $ln$. In the finite field case, $lA$ is always a subset of
the set of residues modulo $n$ and so its cardinality is always at
most $n$, no matter how large $l$ is. This suggests that we should
gain an extra factor $l$ in the assumption of the theorem and that
has turned out to be  indeed the case. The analogue of Theorem
\ref{111} is as follows

\begin{theorem}  \label{111modn} For any fixed positive integer $d$
there are positive constants $C$ and $c$ depending on $d$
 such that the following holds. Let $n $
be a prime and $l$ be a positive integer and   $A$ be a set of
residues modulo $n$ such that  $l^{d+1} |A| \ge C n$. Then the
sumset $lA$ (modulo $n$) contains an arithmetic progression
(modulo $n$) of length $\min \{n, cl |A|^{1/d} \}$.
\end{theorem}

There are two modifications in Theorem \ref{111modn} (compared
with Theorem \ref{111}). First we changed $l^{d}$ to $l^{d+1}$,
which is consistent with the above discussion.  Second, we changed
the lower bound from $ cl |A|^{1/d}$ to $\min \{ n, cl |A|^{1/d}
\}$. This modification is natural and justified, as $lA$ can have
at most $n$ elements. We shall comment on this at the end of  the
next paragraph.

\vskip2mm

The proof Theorem \ref{111modn} is the same as the proof of
Theorem \ref{111}, the only place  one needs a (formal)
modification is (\ref{ln}). In this inequality, the rightmost
formula should be $Cn$ instead of $Cln$, which is consistent with
the discussion in the  paragraph preceding Theorem \ref{111modn}.
Freiman's theorem and all lemmas used for the proof of Theorem
\ref{111} hold for residue classes (see \cite{szemvu} for exact
statements). To explain the change in the lower bound, notice that
in the proof of Theorem \ref{111} we actually showed that either
$lA=[ln]$ or $lA$  contains an arithmetic progression of length
$cl|A|^{1/d}$. Its finite field analogue says that either $lA$
contains all residues modulo $n$ or it contains  an arithmetic
progression of length $cl|A|^{1/d}$. In Theorem \ref{111}, it is
unnecessary to state the lower bound as  $\min \{ln, cl|A|^{1/d}
\}$ because $ln$ is always larger than $cl|A|^{1/d}$. On the other
hand, in the finite field case, it makes sense to write $\min \{n,
cl |A|^{1/d} \}$ since  $n$ can be smaller than $cl |A|^{1/d} \}$.

\vskip2mm Theorem \ref{111modn} demonstrates the flexibility of
our method. It is not clear, for instance, how to prove a finite
field version of Theorem \ref{FS} (which is a special case of
Theorem \ref{111})  using the original approaches of Freiman and
S\'ark\"ozi.

\vskip2mm Similar to Theorem \ref{111}, Theorem \ref{111modn} is
sharp. One can modify the general construction in Section 3 to
match the lower bound. This construction also mirrors the extra
term $l$.

\vskip2mm

{\it \noindent A construction modulo $n$.} We present a
modification of the principal construction in Section 3. Now set
$a= \lfloor \frac{(1-\delta/3)n/l}{d |A|^{1/d}} \rfloor$ (notice
the extra $l$ in the nominator) and $b = \lfloor(\frac{n}{dl
|A|^{1/d}})^{1/(d-1)} \rfloor$. Notice that under the assumption
of Theorem \ref{111modn}, (\ref{const1}) stills hold with the new
definition of $a$. We again have two cases:

\vskip2mm

(I) $\sum_{i=1}^d r_i =0 (\mod \,\, \,\, n)$. By the definition of
the $a_i$'s, it follows that $\sum_{i=1}^d r_i b_i=0 (\mod \,\,
\,\, n)$ and $d$ should be at least 3. By the definition of the
$b_i$'s, it follows immediately that

\begin{equation}  \label{lower3} \max_{ 1\le i \le d} |r_i| \ge\min (
\min_{1\le j \le d} \frac{b_j} {\sum_{i=1}^{j-1} b_i },
\frac{n}{\sum_{j=1}^d b_j} ) \ge \frac{1}{2} a^{1/(d-1)} \ge
2l|A|^{1/d}, \end{equation}

\noindent where the last inequality is from (\ref{const1}).

\vskip2mm (II) $\sum_{i=1}^d r_i \neq 0 (\mod \,\, \,\, n)$. In
this case, we have

$$ \sum_{j=1}^d r_j a +\sum_{j=1}^d r_j b_j = pn
$$

\noindent for some integer $p$. If $p=0$, then

\begin{equation}  \label{lower4} \max_{1\le i \le d} |r_i| \ge  \frac{a}
{\sum_{i=1}^{d} b_i} \ge \frac{1}{2} a^{1/(d-1)} \ge 2l|A|^{1/d}.
\end{equation}

\noindent If $p \neq 0$, then

\begin{equation}  \label{lower5} \max_{1\le i \le d} |r_i| \ge \frac{n}{da
+\sum_{j=1}^d b_j} \ge \frac{n}{(d+1)a} \ge  \frac{1}{2}
a^{1/(d-1)} \ge 2l|A|^{1/d}. \end{equation}

\vskip2mm

Without any further explanation, we now state the analogues of
Theorems \ref{211}, \ref{311}  and  \ref{411}.

\begin{theorem}  \label{211modn} For any fixed positive integer $d$
there are positive constants $C$ and $c$ depending on $d$
 such that the following holds. Let $A_1, \dots, A_l$ be sets of residue classes modulo $n$
  of size $|A|$ where
 $l$ and $|A|$ satisfy $l^{d+1} |A| \ge Cn$.
Then $A_1 +\dots +A_l$ either contains all residue classes modulo
$n$ or contains   a proper GAP of rank $d'$ and volume at least
$cl^{d'} |A|$, for some integer $1 \le d' \le d$.
\end{theorem}

\begin{theorem}  \label{311modn} For any fixed positive integer $d$
there are positive constants $C$ and $c$ depending on $d$
 such that the following holds. Let $n $
be a prime and $l$ be a positive integer and   $A$ be a set of
residues modulo $n$ such that  $l^{d+1} |A| \ge C n$. Then  $lA$
either contains all residue classes modulo $n$ or contains a
proper GAP of rank $d'$ and volume at least $cl^{d'} |A|$, for
some integer $1 \le d' \le d$.
\end{theorem}

\begin{theorem}  \label{411modn} For any fixed positive integer $d$
there are positive constants $C$ and $c$ depending on $d$
 such that the following holds. Let $n $
be a prime and $l$ be a positive integer and   $A_1, \dots, A_l$
be sets of residues modulo $n$ such that  $|A_1|= \dots
=|A_l|=|A|$ and $l^{d+1} |A| \ge C n$. Then $A_1
\stackrel{\ast}{+} \dots \stackrel{\ast}{+} A_l$ either contains
all residue classes modulo $n$ or contains a proper GAP of rank
$d'$ and volume at least $cl^{d'} |A|$, for some integer $1 \le d'
\le d$.
\end{theorem}

\subsection {An application}

A set $A$ of residues modulo $n$ is called {\it zero-sum-free} if
none of the subset of $A$ adds up to  zero modulo $n$.
Zero-sum-free sets are objects of considerable interest in
additive number theory (see  Section C of \cite{Guy} and the
references therein). Here we address the following basic question:

\vskip2mm \centerline {\it How many zero-sum-free sets are there
?}

\vskip2mm

We denote by $S_A$ the collection of partial sums of $A$, so $A$
is zero-sum-free if and only if $0 \notin A$. Szemer\'edi
\cite{Sze0} and Olson  \cite{Ols}, answering a question of
Erd\"os, proved that a zero-sum-free set has at most $2n^{1/2}$
elements. This implies that  the number of zero-sum-free sets is
at most

$$\sum_{i=1}^{ \lfloor 2n^{1/2} \rfloor} {n \choose i} =2^{ \Omega
(n^{1/2} \log_2 n)}. $$

It is not hard to give a lower bound of $2^{\Omega (\sqrt {n}})$;
notice that every subset of the interval $[ \lfloor \sqrt{ 2n}-1
\rfloor]$ is zero sum free, since

$$1 +2+\dots + [ \lfloor \sqrt{ 2n}-1
\rfloor] < n. $$

\noindent The number of subsets of the above interval is clearly
 $2^{\Omega (\sqrt {n}})$.

 \vskip2mm In an earlier paper \cite{szemvu}, we succeeded to
 establish a sharp bound, using a weaker version of Theorem
 \ref{411modn}. (To be more precise, what we actually used was a weaker
 version of the finite field analogue of Theorem \ref{1}.)

\begin{theorem}  \label{0sumfree} Let $n$ be a prime.
The number of zero-sum-free sets ($\mod \,\, n$) is

$$2^{(\sqrt {\frac{1}{3}} \pi \log_2 e +o(1)) \sqrt n}.$$ \end{theorem}

This surprising   estimate  might deserve an explanation. To
reveals its origin, let us give a short proof for the lower bound.
We call a set $A$ of positive integers $n$-{\it small} if the sum
of the elements in $A$ is less than $n$.  It is trivial that an
$n$-small set is zero-sum-free. On the other hand, the number of
$n$ small sets is  $2^{(\sqrt {\frac{1}{3}} \pi \log_2 e +o(1))
\sqrt n}$ due to the following lemma, which is a well-known result
in the theory of partitions (see, for instance, Theorem 6.7 in
\cite{And}).

\begin{lemma} \label{partitions} The number of representations of
$n$ as a sum of different positive integers is $2^{(\sqrt
{\frac{1}{3}} \pi \log_2 e +o(1)) \sqrt n}$. Consequently, the
number of $n$-small sets is $$2^{(\sqrt {\frac{1}{3}} \pi \log_2 e
+o(1)) \sqrt n}. $$ \end{lemma}

The hard part of Theorem \ref{0sumfree} is the upper bound. Using
our results on long arithmetic progressions (modulo $n$)  we
managed to show that if $A$ is zero-sum-free and has relatively
many elements (the number of sets with at most $n^{1/2}/ \log_2^2
n$ elements is $2^{o (n^{1/2})}$ so we can ignore these sets),
then $A$ is close to be $n$-small (for the exact statement please
see \cite{szemvu}). The general idea is as follows. Let $A'$ be a
relatively small subset of $A$; our results show that $S_{A'}$
contains a quite long arithmetic progression. We next make many
translations of this arithmetic progression by adding to it
elements from $A\backslash A'$.  If all these translations avoid
$0$, then we have a good chance to deduce a structural property of
$A$ and it turned out that typically $A$ should look like a
$n$-small set. A similar argument can be applied to determine  the
number of $x$-sum-free sets, for any non-zero residue class $x$.
Trying not to spoil the fun, we do not state the theorem here (it
can be found in \cite{szemvu}), but let us mention that the bound
for non-zero $x$ is different from the bound in Theorem
\ref{0sumfree}. Guessing this bound is a nice puzzle the reader
who bears with us until this very end might enjoy.

\end{document}